\documentclass[dvipsnames]{amsart}
\pdfoutput=1
\usepackage[left=1.15in,right=1.15in,top=1in]{geometry}

\makeatletter
\renewcommand{\email}[2][]{%
  \ifx\emails\@empty\relax\else{\g@addto@macro\emails{,\space}}\fi%
  \@ifnotempty{#1}{\g@addto@macro\emails{\textrm{(#1)}\space}}%
  \g@addto@macro\emails{#2}%
}
\makeatother

\usepackage{enumerate}
\usepackage{amsmath,amssymb,graphicx,bbm}
\usepackage{verbatim,rotating}
\usepackage{epstopdf}
\usepackage{mathrsfs,fancyhdr,xcolor}

\usepackage[ruled,vlined]{algorithm2e}
\usepackage{subfig}

\DeclareMathOperator*{\argmin}{argmin}

\def\R{\mathbb{R}}
\def\N{\mathbb{N}}

\newcommand{\dx}[1]{\mathrm{d} #1}
\usepackage{mathtools}
\newcommand{\bs}[1]{\boldsymbol{#1}}

\usepackage[foot]{amsaddr}

\title[Induced sampling for data-driven polynomial Chaos]
{
{Sparse approximation of data-driven Polynomial Chaos expansions: an induced sampling approach}
}
\thanks{L.~Guo is
supported by the NSF of China (No.11671265). This material is based upon work supported by the National Science Foundation
under Grant No. DMS-1439786 and by the Simons Foundation Grant No. 50736 while
A. Narayan was in residence at the Institute for Computational and Experimental
Research in Mathematics in Providence, RI, during the "Model and dimension
reduction in uncertain and dynamic systems" program. T.~Zhou is partially supported by the NSF of China (under grant numbers 11688101, 11571351, and 11731006), science challenge project (No. TZ2018001), and the youth innovation promotion association (CAS).}

\author{Ling Guo}
\author{Akil Narayan}
\author{Yongle Liu}
\author{Tao Zhou}

\address[L.~Guo]{Department of Mathematics, Shanghai Normal University, Shanghai, China. Email: lguo@shnu.edu.cn.}

\address[A.~Narayan]{Department of Mathematics, and Scientific Computing and Imaging Institute, University of Utah, Salt Lake City, UT. Email: akil@sci.utah.edu.}
\address[Y.~Liu]{Department of Mathematics, Southern University of Science and Technology, Shenzhen, China. Email: 11749318@mail.sustc.edu.cn.}
\address[T.~Zhou]{LSEC, Institute of Computational Mathematics and Scientific/Engineering
Computing, Academy of Mathematics and Systems Science, Chinese Academy of Sciences, Beijing, China. Email: tzhou@lsec.cc.ac.cn.}

\begin{document}
\maketitle

\graphicspath{{figures/}}

\begin{abstract}
  One of the open problems in the field of forward uncertainty quantification (UQ) is the ability to form accurate assessments of uncertainty having only incomplete information about the distribution of random inputs. Another challenge is to efficiently make use of limited training data for UQ predictions of complex engineering problems, particularly with high dimensional random parameters. We address these challenges by combining data-driven polynomial chaos expansions with a recently developed preconditioned sparse approximation approach for UQ problems. The first task in this two-step process is to employ the procedure developed in (Ahlfeld et al. 2016, \cite{Ahlfeld_2016SAMBA}) to construct an ``arbitrary" polynomial chaos expansion basis using a finite number of statistical moments of the random inputs. The second step is a novel procedure to effect sparse approximation via $\ell^1$ minimization in order to quantify the forward uncertainty. To enhance the performance of the preconditioned $\ell^1$ minimization problem, we sample from the so-called induced distribution, instead of using Monte Carlo (MC) sampling from the original, unknown probability measure. We demonstrate on test problems that induced sampling is a competitive and often better choice compared with sampling from asymptotically optimal measures (such as the equilibrium measure) when we have incomplete information about the distribution. We demonstrate the capacity of the proposed induced sampling algorithm via sparse representation with limited data on test functions, and on a Kirchoff plating bending problem with random Young's modulus.
\end{abstract}



\pagestyle{myheadings}
\thispagestyle{plain}

\section{Introduction}
The main purpose of uncertainty quantification (UQ) is to quantify the effect of various sources of randomness on output model predictions. An effective and popular approach is to
build an approximation for the map between the random input and the model output,i.e., the quantities of interest (Q$\circ$I). Many techniques have been developed to construct such approximations. When the distribution of the random input variables are known, the generalized Polynomial Chaos (gPC) \cite{Wiener_1938homochaos,Xiu_2002Wiener} expansion is a popular approach. The basic idea is to build a polynomial approximation of the QoI and then the goal is to compute the polynomial expansion coefficients. There are intrusive and non-intrusive approaches to compute the unknown gPC coefficients. In this paper, we focus on the non-intrusive stochastic collocation method, which constructs a global polynomial approximations after running the model with a sample set of the random variables. We refer to \cite{Narayan_2015SCUMM,Tang_2015redeUQ} and references therein for recent developments for stochastic collocation methods.

However, in real-world scientific applications, one of the main challenges for UQ is the lack of explicit knowledge of the random inputs. Recently, extension of the gPC method to more general input distributions or arbitrary distributions have been investigated. These methods include Multi-element generalized polynomial chaos (ME-gPC) using local polynomial expansion \cite{Prempraneerach_2010Uq,Wan_2008MEgPC}, Global polynomial expansions based on Gram-Schmidt orthogonalization \cite{Witteveen_07uncerapc,Witteveen_2006MauGS}, and moment match method to deal with arbitrary distributions (termed aPC) \cite{Oladyskin_2012Datadriven}. There has also been extensive work on constructing polynomial approximation for dependent random variables or arbitrary distributions, we refer to the recent paper \cite{Jakeman2019} and references therein.

In this study, we focus on the aPC approach for the situation when one has incomplete input information, especially when only sample locations are given. It has been shown by Oladyshkin and Nowak in \cite{Oladyshkin_2011concept} that aPC approach can propagate the moments and thus it offers reliable results with limited input data. The procedure is straightforward to implement. However, it will become ill-conditioned when the polynomial order is large. In \cite{Ahlfeld_2016SAMBA}, the authors proposed an algorithm that can reduce the amount of ill-conditioning (yet can still be ill-conditioned) by calculating the required quantities directly using only matrix operations performed on the Hankel matrix of moments. A sparse grid approach based on the Smolyak's algorithm was proposed to solve the expansion coefficients.

Another challenge for performing UQ predictions is the limited size of available training data since a single simulation for complex systems may require a large amount of computational resources and can be very time consuming. In such basis, the number of available training samples will be much smaller than the complexity of the desired predictor. Recently, based on the idea of compressive sensing \cite{Candes_2006nos, Candes_2006Stablesrec, Donoho_2006Cs, Donoho_2006srs}, stochastic collocation methods via $\ell^1$ minimization \cite{Doostan_2011nonadapted,Yan_2012Sc,Rauhut_2012sparseLegen, mathelin_gallivan_2012} have been shown to be efficient for computing sparse approximations to gPC coefficients from a small amount of training data. However, the success of this approach depends on the strategy used to generate the training data. A popular approach is a ``standard" Monte Carlo approach that would draw samples from the probability density $\omega$ of the random input variables. To further improve the accuracy of the function approximation with high-degree polynomial approximation, numerous sampling strategies have been proposed recently. In \cite{Rauhut_2012sparseLegen}, the authors used Chebyshev sampling for sparse Legendre polynomial approximations. The high-dimensional case was studied in \cite{Yan_2012Sc} and the theoretical results show that the accuracy degrades with increasing parameter dimension. To cope with this limitation, a coherence optimal sampling strategy is given by the authors in \cite{Hampton_2015Cs}.

A general sampling strategy called Christoffel Sparse Approximation (CSA) is developed in \cite{Jakeman_2016generalizedsample}. The sampling algorithm performs well and can be applied for both bounded and unbounded random inputs. The essential idea of the CSA algorithm is to sample from a certain biased distribution (which is an \emph{equilibrium measure}) that is associated with the random parameter density and state space. The success of the CSA algorithm relies on the idea that sampling measures from which stable least squares estimators can be built are also useful for performing sparse approximation. However, an explicit form for the equilibrium measure is not known in general; alternatively the \emph{induced distribution} introduced in \cite{cohen_optimal_2017} is an attractive sampler for its optimal stability properties. Algorithms for generating samples from fairly general classes of induced distributions are also available \cite{IOPD}. However, such algorithms are less helpful when the underlying distribution is not known.

In cases of unknown distributions, it is common to develop a sampling strategy in a \emph{data-driven} method, meaning that a sampling method is devised using a finite number of observations of the random variable. Recently, a data driven polynomial expansion with weighted least-square approach was proposed in \cite{GLZ_ddpcls} for UQ problems with arbitrary random inputs. However, studies in this direction are at present relatively limited in UQ settings.

The method proposed in this paper combines the idea of data-driven aPC from \cite{Oladyskin_2012Datadriven} with sparse polynomial approximation. The algorithm is novel since the available sparse approximation sampling methods require knowledge of the exact distribution, which we do not have. Inspired by the computation of induced polynomial distributions, the purpose of this paper is to provide an efficient induced sampling strategy for sparse approximation of certain types of arbitrary polynomial chaos expansions. The aPC approach with induced sampling proceeds in two steps: we first construct the aPC basis with a moment-matching method; then a Christoffel-function preconditioned $\ell^1$-minimization method is used to perform sparse approximation seeking the dominant aPC expansion coefficients. A sampling strategy that generates samples from the induced measure is employed to guarantee the accuracy of the sparse approximation when limited data is available. Compared with an equilibrium measure sampling strategy, our induced distribution sampling strategy will not generate samples out of the support of the discrete measure defined by the training samples. The numerical results show the superiority of the induced sampling strategy over CSA and oracle Monte Carlo sampling for several test problems.

The rest of this paper is organized as follows. In Section 2, we introduce the traditional gPC approach and the $\ell^1$ minimization method. The moment match data driven basis is given in Section 3 followed by the induced sampling strategy.  Numerical experiments are then shown in Section 4 to indicate the applicable and effectiveness of our approach. Finally, we give some concluding remarks in Section 5.

\section{Background}

In this section, we first briefly review the generalized polynomial chaos (gPC) expansion approach to quantify uncertainties in mathematical models with independent random inputs. Then we give a few facts about the theory of compressive sensing, which has been well-studied and proved to be an efficient method to find sparse solutions to underdetermined linear system \cite{Candes_2008Rip,Davies_2010RICLp,Donoho_2006srs}.

\subsection{Generalized polynomial chaos}\label{ssec:gpc}

We assume $Z=(Z_1,Z_2,\ldots,Z_d)$ is a vector of independent random input parameters in the model and $f(Z): \mathbb{R}^d\rightarrow \mathbb{R}$ is the model output (prediction).
For each $Z_i$ in $\Gamma_i \subset \R$ it admits a marginal probability density $\rho_i$. Under the assumption of independence, we have the joint density function for $Z$ denoted by $\omega(z)= \prod_{i=1}^d \omega_i(z_i): \Gamma\rightarrow \mathbb{R}^+$ with $\Gamma:= \prod_{i=1}^d\Gamma_i \subset \mathbb{R}^d$.
The gPC approach aims to construct a polynomial approximation of the model output $f(z)$ as follows:

\begin{equation}\label{eq:multi-expan}
f(z)\approx f_N(z)=\sum_{\boldsymbol{\lambda}\in\Lambda}c_{\boldsymbol{\lambda}}\Phi_{\boldsymbol{\lambda}}(z),
\end{equation}
where $\boldsymbol{\lambda}=\{{\lambda}_1,{\lambda}_2,\ldots,{\lambda}_d\} \in \N_0^d$ is a multi-index, $\Lambda$ is a multi-index set (i.e., a subset of $\N_0^d$) that has size $|\Lambda| = N$, and $\Phi_{\boldsymbol{\lambda}}$ are basis functions that are typically chosen to be orthonormal with respect to the density $\omega(z),$, i.e.,
\begin{equation}
(\Phi_{\boldsymbol{\lambda}}(z),\Phi_{\boldsymbol{\theta}}(z))_{\omega} \coloneqq \int_\Gamma \!\!  \Phi_{\boldsymbol{\lambda}}(z) \Phi_{\boldsymbol{\theta}}(z)\omega(z)dz = \delta_{\boldsymbol{\lambda},\boldsymbol{\theta}}, \quad \boldsymbol{\lambda},\boldsymbol{\theta}\in\Lambda.
\end{equation}
While any basis for \eqref{eq:multi-expan} spanning the same space is acceptable, the orthonormality assumption is made for two main reasons: (i) statistics of the approximation are efficient and straightforward to compute from coefficients of an orthonormal expansion, and (ii) computational algorithms that yield the gPC coefficients for an orthonormal family are frequently more numerically stable than those that use a non-orthogonal family. The assumption of independence of the components $Z_i$ allows us to choose $\Phi_{\boldsymbol{\lambda}}$ as tensor-products of the univariate orthogonal polynomials in each direction, i.e.,
\begin{equation}\label{eq:d-variate-polynomial}
\Phi_{\boldsymbol\lambda}(z)=\prod_{i=1}^d\phi_{{\lambda}_i}^i(z_i) \quad  \textmd{with} \quad
  \int_{\Gamma_i}\!\! \phi^i_{{\lambda}_k}(z_i) \phi^i_{{\lambda}_l}(z_i) \omega_i(z_i) dz_i = \delta_{k,l}.
\end{equation}
The main advantage of reducing this problem to products of one-dimensional orthonormal polynomials is that there are stable and efficient procedures to generate these polynomials in one dimension, namely the three-term recurrence relation. For dimension $i$ there are coefficients $\{a^i_j\}_{j \geq 1} \subset \R$ and $\{b^i_j\}_{j \geq 0} \subset (0, \infty)$ that depend only on $\omega_i$ such that
\begin{align}\label{eq:ttr}
  z_i \phi^i_{\ell}(z_i) &= b^i_{\ell+1} \phi^i_{\ell+1}(z_i) + a^i_{\ell+1} \phi^i(z_i) + b^i_{\ell} \phi^i_{\ell-1}(z_i), & \ell &\geq 0,
\end{align}
with $\phi^i_0(z_i) \coloneqq 1/\sqrt{b^i_0}$. Knowledge of these coefficients allows one to stably evaluate and manipulate univariate orthogonal polynomials.

In practice, $\Lambda$ is a finite multi-index set, and is particular chosen as a monotone or downward-closed set, i.e.,
\begin{align*}
  \bs{\lambda} \in \Lambda  \Rightarrow \bs{\lambda} - \bs{e}_j \in \Lambda, \hskip 3pt j = 1, \ldots, d,
\end{align*}
where $\bs{e}_j \in \R^d$ is a vector whose only nonzero entry is 1 in position $j$. Several types of canonical index sets are monotone, including the total degree space,
\begin{equation*}
\Lambda = \Lambda_k^{\textrm{TD}} \coloneqq \left\{\boldsymbol\lambda \;\; \big| \;\;  |\boldsymbol\lambda|_1=\sum_{i=1}^d \lambda \leq k \right\}.
\end{equation*}
In this paper we will focus our examples on this particular family of index sets, but our methodology does not require $\Lambda$ to be total-degree, nor does it require $\Lambda$ to be downward-closed.
The associated finite-dimensional subspace of functions to which $f_N$ belongs is the total degree polynomial space,
\begin{align}\label{eq:TD}
P(\Lambda) \coloneqq \mathrm{span}\left\{ \Phi_{\boldsymbol\lambda} \;\; \big| \;\;  \boldsymbol\lambda \in \Lambda_k^{\textrm{TD}}\right\}.
\end{align}
The dimension of this total degree space is $N=\bigg (\begin{array}{c}
                                                       d+k \\
                                                       d
                                                     \end{array}\bigg )
$. For convenience, we can place an appropriate order on the multi-indices, i.e.,
\begin{align}\label{eq:index-ordering}
  \left\{ \boldsymbol{\lambda} \;\; | \;\; \boldsymbol{\lambda} \in \Lambda \right\} \longleftrightarrow \left\{ 1, \ldots, N \right\}.
\end{align}
Thus we can instead identify the basis functions using a scalar index,
\begin{equation}\label{eq:multi-index-single}
\{\Phi_{\boldsymbol{\lambda}}(\xi)\}_{\boldsymbol{\lambda}\in\Lambda}\Leftrightarrow \{\Phi_j(\xi)\}_{j=1}^{N}.
\end{equation}
Hereafter, we will use the single index $\{j=1,2,\ldots,N\}$ instead of the multi-indices. Therefore, the gPC approximation (\ref{eq:multi-expan}) can be written as
\begin{equation}\label{eq:finite-N-expan}
f(z)\approx f_N(z)=\sum_{j=1}^Nc_j\Phi_j(z).
\end{equation}
Many efficient numerical techniques have been developed to estimate the expansion coefficients $\{c_j\}_{j=1}^N$, such as the intrusive stochastic Galerkin methods \cite{Xiu_2002Wiener} and the non-intrusive collocation methods \cite{Eldred_2009nonintrusive,XiuH}. One practical restriction of these approaches is that amount of data or observations of $f$ required to construct the gPC emulator $f_N$ is at least $N = |\Lambda|$. Typically $N$ can be very large, and for example, our total degree set $\Lambda^{\mathrm{TD}}_k$ has cardinality growing nearly exponentially with the dimension $d$, and so the previously mentioned techniques would require a large amount of data. An alternative class of strategies that has gained attention in recent years leverages sparsity in the expansion coefficients $c_j$ that frequently surfaces in applications. Such sparse approximation strategies often utilize $\ell^1$ minimization algorithms coming from the compressed sensing community; these algorithms can construct an accurate $f_N$ with a number of observations that scale only with the number of nonzero expansion coefficients.

\subsection{Sparse approximation via $\ell^1$ minimization}\label{ssec:sparse}
The model $f_N$ in \eqref{eq:finite-N-expan} can have a large number of unknown coefficients $c_j$ to compute, yet the underlying system $f$ is frequently computationally expensive and we can afford only a limited number of simulations. In such cases we do not have enough available data to approximate every coefficient, so instead we hope to recover the dominant values. To make this precise, suppose we are given a set of $M$ unstructured realizations $\{z^{(j)}\}_{j=1}^M$ and the corresponding outputs $b=(f(z^{(1)}),...f(z^{(M)}))^{T}$, we seek a solution vector $c$ to the linear problem,
\begin{equation*}\label{eq:gPC_L0}
Ac\approx b,
\end{equation*}
where $c=(c_1,\ldots, c_N)^T \in \mathbb{R}^N$ is the coefficient vector in \eqref{eq:finite-N-expan}, and $A \in \mathbb{R}^{M \times N}$ denotes the measurement matrix, which is written as
 \begin{eqnarray}\label{eq:matrixelem}
   A=(a_{ij})_{1\leq i\leq M,1\leq j\leq N}, \quad a_{ij}=\Phi_{j}(z^{(i)}).
 \end{eqnarray}
In the case when our unknown coefficients (vastly) outnumber the available data, we have $M \ll N$ and so the linear system $Ac\approx b$ is underdetermined, so we can only recover partial information about $c$. One strategy to recover this information is through a sparsity assumption on $c$ and using computational tools for solving the underdetermined problem in compressed sensing. In particular, we can use the well-studied $\ell^1$ minimization framework to compute the gPC coefficients by solving the optimization problem
\begin{equation}\label{eq:gPC_L1}
\min \|c\|_1 \quad \text{subject to}\ Ac=b.
\end{equation}
This $\l_1$ minimization problem is also referred to as a basis pursuit optimization.  When the data $b$ is contaminated by noise, the constraint in \eqref{eq:gPC_L1} is relaxed to obtain the basis pursuit denoising problem. In this paper, we focus on the noiseless interpolation-type constraints, i.e the basis pursuit problem, which can be solved with efficient algorithms from convex optimization.

\section{Method}

The overall problem we consider is as follows: Suppose we are given $Q$ empirical samples and weights corresponding to the random variable $Z$. Denote these empirical samples as $S = (z^{(1)}, \ldots, z^{(Q)} ) \subset \R^d$, with associated weights $(w_1, \ldots, w_Q)$. For example, if $S$ contains iid samples drawn at random from $\omega$, then assigning the Monte Carlo weights $w_j = 1/Q$ is reasonable. Given this information, we seek to formulate a procedure that specifies $M \ll Q$ sample locations where data is queried for use in the sparse approximation strategy of section \ref{ssec:sparse}. We accomplish this in two steps. Sections \ref{ssec:apc-1d} and \ref{ssec:apc-nd} perform the first step wherein an approximation to the underlying (unknown) distribution is implicitly constructed. The second step proposed in Section \ref{ssec:l1} uses this constructed distribution to specify sample locations and subsequently formulate a preconditioned version of \eqref{eq:gPC_L0} that accomplishes function approximation.

\subsection{Data-driven polynomial chaos: a moment-based approach }\label{ssec:apc-1d}
In this subsection, we will present the arbitrary polynomial chaos (aPC) approach for one-dimensional densities, whose goal is to compute a distribution for the random input given incomplete (moment) information from the true distribution. This strategy was investigated in \cite{Oladyskin_2012Datadriven,Ahlfeld_2016SAMBA} and can handle the situation when only a discrete number of sample locations (or only moment information) is available that characterize the random input. While there are several other types of gPC approaches to cope with such epistemic uncertainty, e.g. \cite{Witteveen_07uncerapc,Witteveen_2006MauGS,Ernst_2012gPCconvergen,soize_2004physicalrandom,Oladyskin_2012Datadriven,Ahlfeld_2016SAMBA}, our focus in this paper is on the aPC approach.

The starting point for the aPC approach is to assume that some (possibly approximate) moments of the unknown distribution are given. In our setting, we assume something slightly different, namely that we are given $Q$ empirical samples of a scalar, real-valued random variable $\xi$. These empirical samples are $\Xi \coloneqq (\xi_1, \ldots, \xi_Q )$ with associated (non-negative) weights $(w_1, \ldots, w_Q)$, which we assume satisfy $\sum_{j=1}^Q w_j = 1$. Our goal is to map the inputs $\Xi$ along with the weights to the orthogonal polynomial three-term recurrence coefficients in \eqref{eq:ttr}. This implicitly produces a set of polynomials $\{\phi_j\}_{j \geq 1}$ in the variable $\xi$ that are orthogonal with respect to a measure that is implicitly defined by the empirical samples. The standard way to address this procedure is as a moment-matching problem.

To begin, we compute un-centered polynomial moments using our empirical samples,
\begin{equation}\label{eq:samples}
\nu_{k}= \sum_{j=1}^Q w_j \xi_j^k, \ \ \ k=0,1,\ldots.
\end{equation}

Recall that our goal is construct an orthogonal polynomial basis, so first for a fixed $K$ the polynomial $\phi_K$ has the expansion
\begin{align}\label{eq:phiK-expansion}
  \phi_K(\xi)=\sum_{j=0}^{K}\beta_j\xi^j, 
\end{align}
We now assume that $\rho$ is some probability density function that has the moments specified in \eqref{eq:samples}, and that $\{\phi_j\}_{j \geq 0}$ is an orthonormal polynomial family with respect to $\rho$. Then multiplying both sides of \eqref{eq:phiK-expansion} by $\xi^\ell$ and taking the $\rho$-expectation for $\ell = 0, \ldots, K$ yields the following linear system for the coefficients $\beta_j$:
\begin{gather}\label{gsshiqi}
\begin{bmatrix}
\nu_{0}&\nu_{1}&\cdots&\nu_{K}\\
\nu_{1}&\nu_{2}&\cdots&\nu_{K+1}\\
\vdots&\vdots&\vdots&\vdots\\
\nu_{K-1}&\nu_{K}&\cdots&\nu_{2K-1}\\
0&0&\cdots&1
\end{bmatrix}
\begin{bmatrix}
\beta_{0}\\
\beta_{1}\\
\vdots\\
\beta_{K-1}\\
\beta_{K}
\end{bmatrix}=\begin{bmatrix}
0\\
0\\
\vdots\\
0\\
1
\end{bmatrix}.
\end{gather}
The polynomial coefficients $\{\beta_j\}_{j=0}^{K}$ can be computed directly by solving this linear system. Note that this also reveals that the ``true" underlying density $\rho$ need only be known up to its first $2 K$ moments in order to compute all polynomials up to degree $K$.
Unfortunately, the moment matrix on the left-hand side is frequently ill-conditioned for large $K$. To mitigate this ill-conditioning, an alternative approach was proposed in \cite{Ahlfeld_2016SAMBA} by considering matrix operations on the Hankel matrix of moments, defined as follows:
\begin{gather}\label{eq:Hankelmatrix}
\mathbf{H}=
\begin{bmatrix}
\nu_{0}&\nu_{1}&\cdots&\nu_{k}\\
\nu_{1}&\nu_{2}&\cdots&\nu_{k+1}\\
\vdots&\vdots&\vdots&\vdots\\
\nu_{K}&\nu_{K+1}&\cdots&\nu_{2K}
\end{bmatrix}.
\end{gather}
We require that the set of $M$ samples is determinate in the Hamburger sense which means that all the corresponding quadratic forms are strictly positive, that is $\text{det}(\mathbf{H})>0$. A sufficient condition for the moments defined by \eqref{eq:samples} to result in a $(K+1) \times (K+1)$ determinate problem is that the empirical set $\Xi$ must contain at least $K+1$ distinct samples.
When the problem is determinate, we can perform a Cholesky decomposition, $\mathbf{H}=\mathbf{R}^{\top}\mathbf{R}$ with
\begin{gather}\label{eq:Cholmatrix}
\mathbf{R}=
\begin{bmatrix}
r_{11}&r_{12}&\cdots&r_{1,K+1}\\
      &r_{22}&\cdots&r_{2,K+1}\\
      &  &\ddots&\vdots\\
      &   &     &r_{K+1,K+1}
\end{bmatrix}.
\end{gather}
Then the entries of the matrix $\mathbf{R}$ can be used to form an orthogonal system of polynomials by the Mysovskih theorem \cite{Mysovskikh_1968cubature}. Moreover, the three-term recurrence coefficients $\{a_j\}_{j \geq 1}$ and $\{b_j\}_{j \geq 0}$ associated to a $\xi$ version of \eqref{eq:ttr} can be computed from the following explicit analytic formulas \cite{Golub_1968Gqr}:
\begin{align}\label{coef}
a_j=\frac{r_{j,j+1}}{r_{j,j}}-\frac{r_{j-1,j}}{r_{j-1,j-1}}, \quad b_{j}=\frac{r_{j+1,j+1}}{r_{j,j}},
\end{align}
where $r_{0,0}=1$ and $r_{0,1}=0$. This procedure, the univariate aPC approach, empirically is much less ill-conditioned compared to directly computing solutions to the system \eqref{gsshiqi}.

Algorithm \ref{alg_aPC} summarizes the procedure of the univariate aPC basis construction via the moment method when the underlying measure $\rho$ is implicitly
represented by a sample set $\Xi$.

\begin{algorithm}
\SetAlgoNoLine
\caption{Moment-based arbitrary polynomial construction}\label{alg_aPC}
\begin{enumerate}
  \item Input: data and corresponding weights, $\{\xi_j,w_j\}_{j=1}^Q$, maximal degree $K$
  \item Compute moments according to equations \eqref{eq:samples};
  \item Use equations \eqref{eq:Hankelmatrix},\eqref{eq:Cholmatrix} and \eqref{coef}, we can compute the polynomial recurrence coefficients $\{a_j,b_j\}$;
  \item (Optional) Evaluate polynomials according to \eqref{eq:ttr}.
\end{enumerate}

\end{algorithm}

\subsection{Tensorization of the univariate aPC approach}\label{ssec:apc-nd}
We now briefly describe a multivariate version of the aPC algorithm from the previous section. Our approach is to implicitly construct a tensorial measure whose marginals match the marginals of the empirical set $S$ for a specified set of moments. We recall our notation for $S$ and its elements:
\begin{align*}
  S &= \left\{ z^{(1)}, \ldots, z^{(Q)} \right\} \subset \R^d,  & z^{(j)} &= \left( z^{(j)}_1, \ldots, z^{(j)}_d \right) \in \R^d.
\end{align*}
With this notation, we introduce the marginalization of the set $S$ in each of the coordinates:
\begin{align*}
  \Xi_i &= \left\{ z^{(1)}_i, \ldots, z^{(Q)}_i \right\} \subset \R, & i &= 1, \ldots, d.
\end{align*}
Associated to each of these sets we assign the weights $w_j = \frac{1}{Q}$ for $j = 1, \ldots, Q$. (The weights do not depend on the dimensional index $i$.) Then given some user-specified maximal degree $K \in \N$ along with $\Xi_i$ and the uniform weights, we use Algorithm \ref{alg_aPC} described in Section \ref{ssec:apc-1d} to construct three-term recurrence coefficients $b_0^i$ and $\left(a_j^i, b_j^i\right)_{j = 1}^K$. In particular, this implicitly defines the polynomials $\{\phi^i_j\}_{j=0}^K$ through the three-term relation \eqref{eq:ttr}.

We repeat this process for each $i = 1, \ldots, d$, resulting in $d$ univariate orthogonal polynomial families that match moments each respective sample marginalization $\Xi_i$. Through \eqref{eq:d-variate-polynomial}, this also defines our tensorized multivariate orthogonal polynomial family. Note that our multivariate polynomials $\Phi_{\bs{\lambda}}$ constructed in this way do \emph{not} respect the moments of $S$. However, they do respect the moments of the tensorized univariate marginalization:
\begin{align*}
  \int_\Gamma \Phi_{\bs{\lambda}}(z) \Phi_{\bs{\theta}}(z) \dx{\nu}(z) &= \delta_{\bs{\lambda}, \bs{\theta}}, & \bs{\lambda}, \bs{\theta} &\in \Lambda^{\mathrm{TD}}_K,
\end{align*}
where $\nu$ is the tensorial measure,
\begin{align}\label{eq:nu-def}
  \dx{\nu}(z) &= \otimes_{i=1}^d \dx{\nu}_i(z_i), & \dx{\nu}_i(z_i) &\coloneqq \frac{1}{Q} \sum_{j=1}^Q \delta_{z^{(j)}_i}(z),
\end{align}
with $\delta_{z_0}(z)$ the univariate Dirac mass centered at $z = z_0$. Thus, our tensorized procedure incurs an additional approximation penalty that is essentially equivalent to ``tensorizing" the set $S$. The construction of the basis functions $\left\{\Phi_{\bs{\lambda}}\right\}_{|\bs{\lambda}| \leq K}$ accomplished above concludes the first step of our data-driven procedure. As described in Section \ref{ssec:gpc}, we will hereafter use an ordered indexing of these basis elements, $\{\Phi_j\}_{j=1}^N$, where $N = |\Lambda_k^{\mathrm{TD}}|$.
The second step of our procedure, described in the next section, performs function approximation using these basis elements.


\subsection{Postprocessing via preconditioned $\ell^1$ minimization}\label{ssec:l1}
Upon completion of the first step of our algorithm at the end of Section \ref{ssec:apc-nd}, we have access to both $d$ univariate sets of orthogonal polynomials $\{\phi_j^i\}_{j=0}^K$ for $i = 1, \ldots, d$. With this, one can perform any standard algorithm that relies on tensorial methods for approximation. For example, a sparse grid method was proposed in \cite{Ahlfeld_2016SAMBA}, where the collocation points are generated based on the data-driven univariate measures $\dx{\nu}_i$. In this section, we shall propose instead to use a preconditioned $\ell^1$ minimization approach.

We propose to solve a modified version of the $\ell^1$ minimization procedure described in Section \ref{ssec:sparse}. Given a function $f$, we generate $M$ subsamples $\{z^{(j)}\}_{j=1}^M$ of the empirical set $S$\footnote{The indexing of samples $z^{(j)}$ in this section do not match the indexing of samples $z^{(j)}$ in Section \ref{ssec:apc-1d} and \ref{ssec:apc-nd}. However, we avoid introducing a new index mapping for notational simplicity, since we believe the meaning here is clear.}. We will describe in this section how this subsampling is performed.

From these $M$ subsamples, we populate the matrix $A$ as in \eqref{eq:matrixelem}, and we construct an $M$-vector $b$ having the entries $b_j = f(z^{(j)})$. We propose to construct an approximation $f_N$ to $f$ as in \eqref{eq:finite-N-expan} where the coefficients $c_j$ are defined by the solution to,
\begin{equation}\label{eq:gPC_preL1}
\min \|c\|_1 \quad \text{subject to}\ \sqrt{W}Ac=\sqrt{W}b,
\end{equation}
where $W=\textmd{diag}(W_1, ... ,W_M)$ is a preconditioning matrix that we will specify. Our selection of the subsamples $z^{(j)}$ and the preconditioning weights $W_j$ are motivated by results in compressed sensing and least squares approximation. We will first describe the details of the approach, and follow this with a motivational discussion.

We define the measure $\mu$ having support $S$ as
\begin{align}\label{eq:mu-def}
  \dx{\mu}(z) &\coloneqq \sum_{j=1}^Q \widetilde{\kappa}_j \delta_{z^{(j)}}, & \widetilde{\kappa}_j &\coloneqq \frac{\kappa(z^{(j)})}{\sum_{q=1}^Q \kappa(z^{(q)})},
\end{align}
with $\kappa(\cdot) = \kappa(\cdot; \Lambda)$ the normalized $\Lambda$-Christoffel function for the tensorial measure $\nu$,
\begin{align}\label{eq:Christoffel}
  \kappa\left(z;\Lambda_K^{\mathrm{TD}}\right) &\coloneqq \frac{1}{N} \sum_{j=1}^N \Phi_j^2(z).
\end{align}
With the measure $\mu$ defined, our subsamples and preconditioning weights for \eqref{eq:gPC_preL1} are computed as follows:
\begin{itemize}
  \item The subsamples $z^{(j)}$ are drawn as $M$ iid samples from $\mu$.
  \item The preconditioning weights are $W_j = 1/\kappa\left(z^{(j)}\right)$.
\end{itemize}

Generating iid samples from the measure $\mu$ is computationally simple since $S$ is a finite set (via, e.g., inverse transform sampling). Note that this provides a complete description of an algorithmic procedure to compute a solution vector $c$ to \eqref{eq:gPC_preL1}, which in turn defines the function approximation $f_N$ in \eqref{eq:finite-N-expan}. This algorithmic procedure is the proposed data-driven aPC approach of this paper.

We now present some discussion and motivation of this approach. For compressed sensing problems that utilize $\ell^1$ minimization, the concept of bounded orthonormal systems can be used to quantify the effectiveness of \eqref{eq:gPC_preL1} as a proxy for a convex $\ell^0$ minimization approach \cite{Rauhut_2010CsSM}. In particular, theory that quantifies error in the ability of \eqref{eq:gPC_preL1} to recover a sparse vector has penalty terms that are proportional to the maximum entry of the matrix $\sqrt{W} A$. Therefore, small maximum magnitudes, ensuring mass is ``spread" evening across the elements of $\sqrt{W} A$, produces both better theory and practical results. For example, this idea was explored in \cite{Rauhut_2012sparseLegen} for sparse univariate Legendre polynomial expansions. With this in mind, one reason why the choice of $W$ is reasonable is that it ensures that all entries of $\sqrt{W} A$ are uniformly bounded, regardless of the sample set $S$ or the polynomial functions $\Phi_j$.

The preconditioning by $\sqrt{W}$ introduces what is effectively a biasing of the measure. In order to retain unbiasedness with respect to the uniform measure on $S$, we must sample from the appropriate modified measure, which is $\mu$. Thus, our aPC approach for \eqref{eq:gPC_preL1} defines sampling and preconditioning based on methods that are known to promote effectively recovery of sparse signals in compressed sensing.

\subsubsection{Relationship to existing methods}\label{sssec:apc-survey}
We devote this section to a brief discussion of existing approaches in the literature that are related to our data-driven aPC approach.

For overdetermined least squares problems, the authors in \cite{GLZ_ddpcls} use data-driven polynomial chaos expansions with weighted least-squares approach to solve uncertainty quantification (UQ) problems. The authors in \cite{Jakeman_2016generalizedsample} utilize precisely the same preconditioning matrix $W$ to solve $\ell^1$ optimization problems, but choose to sample not from $\mu$, but instead from a $K$-asymptotic version of $\mu$, called the pluripotential equilibrium measure. Such a sampling method has also been investigated for overdetermined least-squares problems \cite{Narayan_2017CSALS}. On the other hand, using a biased measure defined by multiplying an existing underlying measure by $\kappa$ (which is how $\mu$ is defined) has become popular in recent years, and was first investigated for least-squares problems in \cite{hampton_coherence_2015,cohen_optimal_2017}. A measure constructed in this way is called an \emph{induced} measure, and computational algorithms for sampling from several multivariate induced distributions exist \cite{IOPD}, and sampling from a multivariate induced measure for gPC approximations is considered in \cite{Jakeman2019} using a different strategy to construct the basis $\Phi_j$. Notice that in all above mentioned references, the input density is known in prior.

With the previous literature survey in mind, our aPC approach can be summarized as follows: We first use the marginalization of $S$ to compute tensorial orthogonal polynomials as basis functions. We then use these functions to define an induced measure $\mu$ on $S$. We perform compressed sensing for sparse approximation by sampling iid from this induced distribution, and retain an unbiased estimate by weighting our samples by the appropriate likelihood function or Radon-Nikodym derivative $1/\kappa$.


\subsubsection{Relationship to the equilibrium measure}
The procedure proposed in \cite{Jakeman_2016generalizedsample} applies to cases when the underlying density $\omega$ is known, but is similar to our proposed aPC approach. The authors in that work note that, for fairly general $\omega$, there is a unique probability measure $\mu_\infty$ such that $\dx{\mu}_\infty = \lim_{\kappa \rightarrow \infty} \omega _\kappa$, where the equality must be interpreted in the appropriate sense. This measure $\mu_\infty$ is called the weighted pluripotential equilibrium measure. While the form of this density is known in some cases, it is not known for general $\omega$ and $\Gamma$. We refer to \cite{Narayan_2017CSALS,Jakeman_2016generalizedsample} for further discussion on explicit formulas and conjectures for this measure. In relevance to this article, the equilibrium measure when $\omega$ has support on the hypercube $[-1,1]^d$ is the product Chebyshev measure,
$$\dx{\mu}_{\infty}(z) = \frac{1}{\pi^d\prod_{i=1}^d \sqrt{1-z_i^2}},$$
which is easily sampled from. In our data-driven aPC approach, the measure $\mu$ we sample from has the advantage that $\mu$ only has support on $S$, so that a valid sample for compressed sensing is guaranteed to be generated. In addition, while the product Chebyshev measure may be optimal in the limit, it is not optimal for a finite degree $K$.

\begin{figure}[ht!]
\centerline{\includegraphics[width=7.0cm]{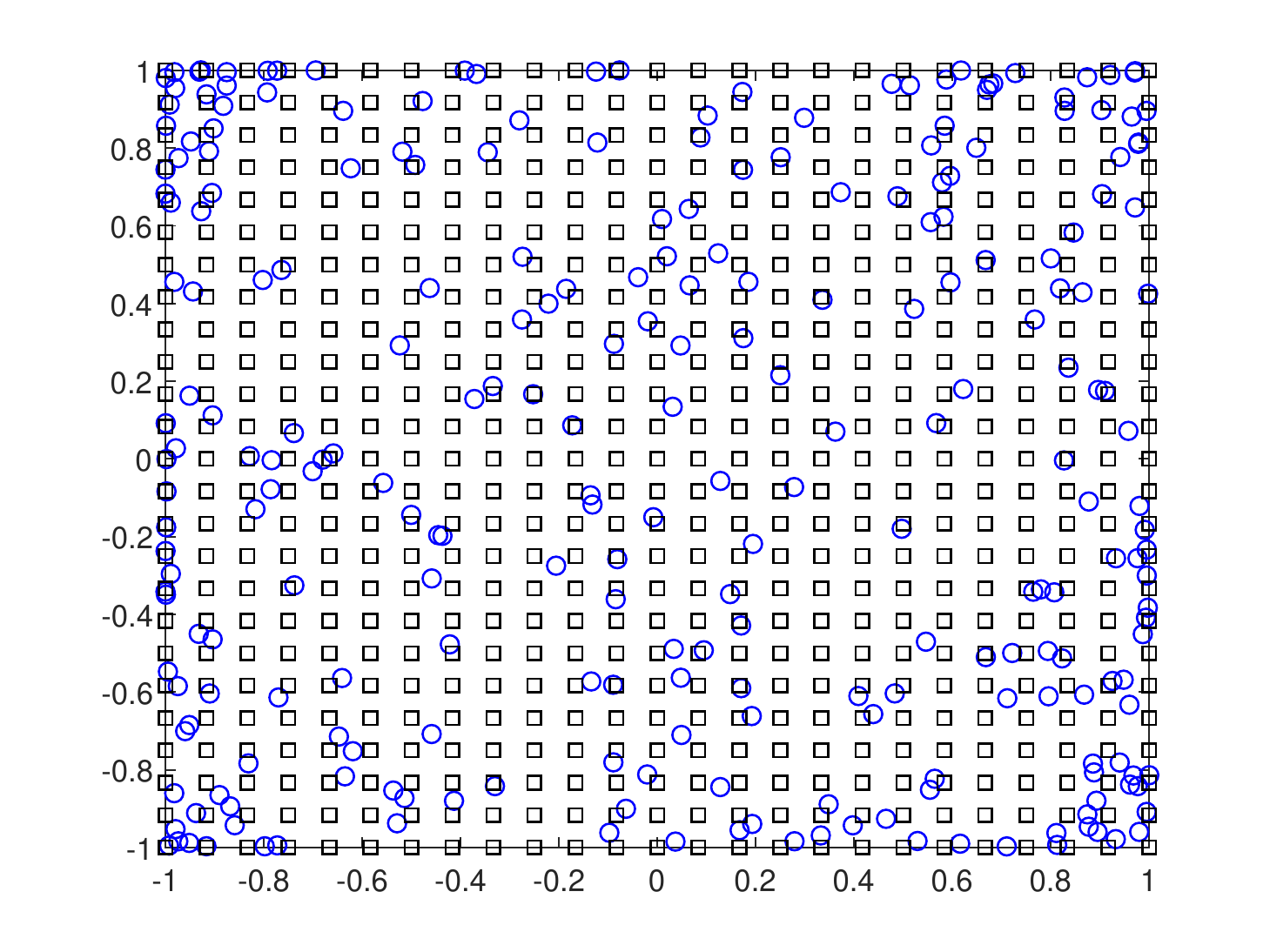}\hspace*{0.5cm}\includegraphics[width=7.0cm]{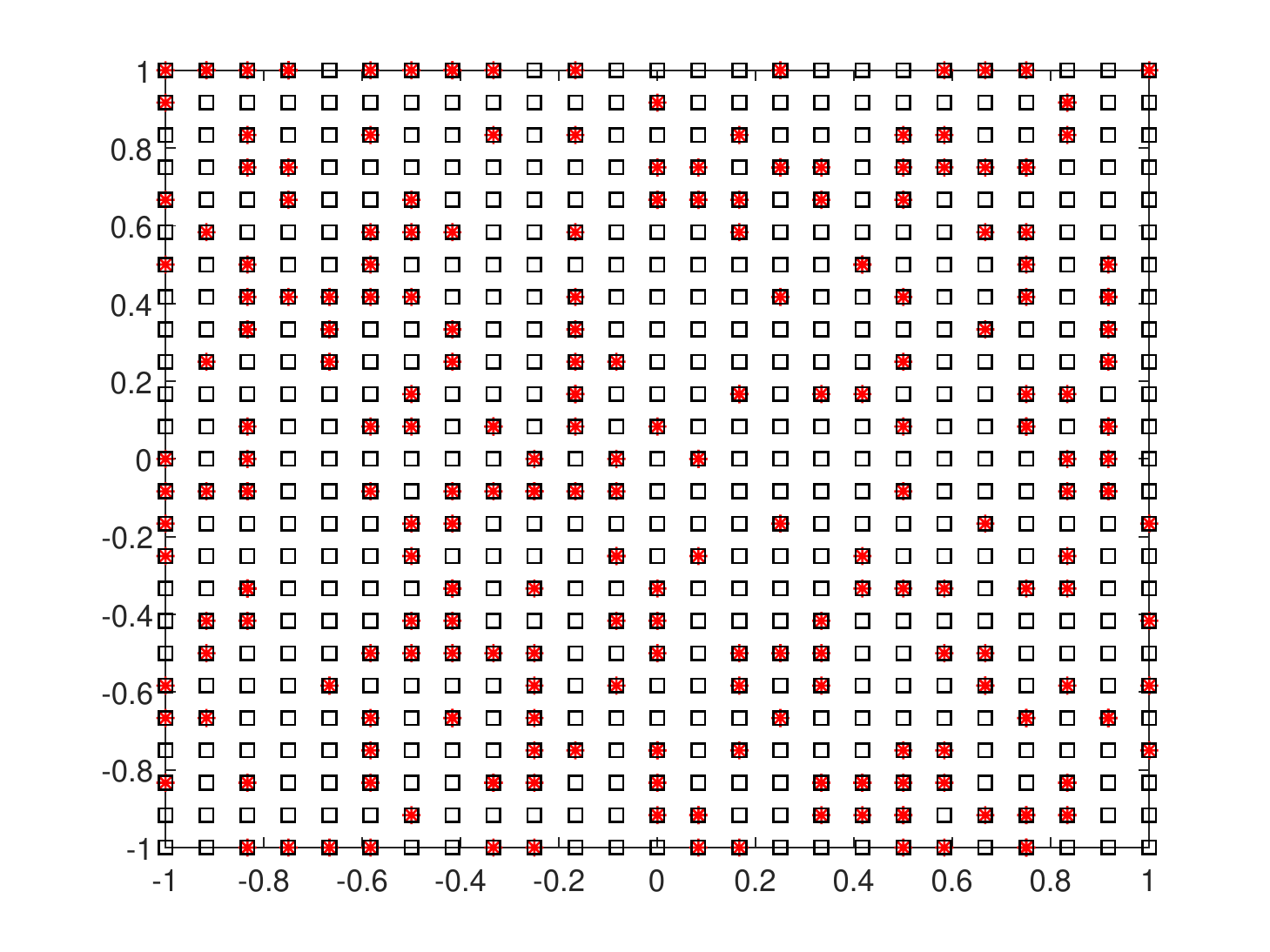}}
  \caption{\sf Sample plots generated via equilibrium measure (left) and induced measure sampling (right). The data set $S$ is tensorial as in \eqref{eq:S-example}, with the $z_1$ samples emulating a Binomial$(24,0.5)$ distibution, and the $z_2$ samples emulating a Poisson($10$) distribution. Left: Samples from the equilibrium (product Chebyshev) measure; Right: samples from the $K=20$ induced sampling measure $\mu$.}\label{sample comparison}
\end{figure}

We visually compare our data-driven aPC approach with equilibrium measure sampling in Figure \ref{sample comparison}. We generate data $S$ as
\begin{align}\label{eq:S-example}
  S = \left\{ \left(z_1^{(j)}, z_2^{(k)}\right)_{1 \leq j, k \leq 24} \right\},
\end{align}
    with associated weights $w_{j,k} = u_j v_k$. The samples $z_1^{(j)}$ are $24$ equispaced points on $[-1,1]$ with weights $u_j$ equal to the probability mass function of a Binomial($24,0.5$) random variable. Thus, the $z_1^{(j)}$ samples emulate samples of a Binomial($24,0.5$) random variable taking values on $[-1,1]$. The samples $z^{(2)}_k$ in turn emulate samples from a truncated Poisson$(10)$ random variable. The samples $z_2^{(k)}$ are 24 equispaced points on $[-1,1]$ with weights $v_k$ equal to the first $24$ values of the probability mass function for a Poisson($10$) random variable. This truncates mass from a Poisson random variable, but the truncated probability is approximately $10^{-5}$. The values $v_k$ are re-normalized so that they sum to $1$. The sample set $S$ is visualized as black squares in Figure \ref{sample comparison}.

The equilibrium measure corresponding to this sample set is just the product Chebyshev measure. In the left pane of Figure \ref{sample comparison} we show the set $S$ versus iid samples generated from the equilibrium measure (the latter visualized as blue dots). The right-hand pane shows iid samples generated from the induced measure $\mu$ for $K=20$, which has support only on $S$, visualized as red dots. Note in particular that the two sampling strategies visually focus effort in different regions. While we expect that the large-$K$ induced measure will behave as the equilibrium measure, this (relatively) small value of $K=20$ is clearly not in this asymptotic regime yet, and so we should expect different approximation results.


Sampling strategies based on the equilibrium measure have only mild dependence on the underlying density $\omega$ and the degree of the polynomial space $K$. In the data-driven framework where the underlining input density is in general unknown, our data-driven aPC sampling approach generates a sampling measure and preconditioning that depends quite explicitly on both $\omega$ and $K$. This suggests that the aPC approach should yield better results compared to equilibrium sampling approaches. We observe this in our numerical results. In addition, the aPC approach is quite flexible: although we only computationally investigate total degree approximation spaces in this article, the data-driven aPC approach applies essentially unchanged to arbitrary finite multi-index sets $\Lambda$.

For more detailed discussions on the induced measure sampling, one can refer to the recent review article \cite{Review}.

\subsection{Data-driven aPC summary}

We summarize the compressed $\ell^1$ approach driven by the induced distribution sampling and arbitrary polynomial chaos construction in Algorithm \ref{alg3}.
\begin{algorithm}
\SetAlgoNoLine
\caption{Induced distribution sampling sparse approximation}\label{alg3}
  Input: the discrete distribution or the data set $\{z^{(j)},w_j\}_{j=1}^Q$, function $f(z)$;
Output: expansion coefficients $c^{\ast}$ such that $f\approx\hat{f}=\sum c^{\ast}_j\Phi_j(z)$
\begin{enumerate}
  \item Construct the arbitrary polynomial space $\{\Phi_j(z)\}$ using Algorithm \ref{alg_aPC} for each dimension;
  \item Generate $M$ iid samples from the induced measure $\mu$
  \item Assemble $b$ with entries $b_i=f(z_i)$ and $A$ with entries $A_{ij}=\Phi_j(z^{(i)})$;
  \item Compute the wieghts $W$ using \eqref{eq:Christoffel};
  \item Compute $c^{\ast}=\argmin_c \|c\|_1$ such that $\sqrt{W}Ac=\sqrt{W}b$.
\end{enumerate}

\end{algorithm}

\section{Results}
In this section, the performance of the proposed data-driven aPC approach, using preconditioned $\ell^1$ minimization with induced distribution sampling, is demonstrated using several examples. In all cases the input to the procedure is a size-$Q$ multivariate sample size $S$ produced as iid samples from a true ``unknown" density $\omega$. The aPC procedure uses the methodology from Sections \ref{ssec:apc-1d}, \ref{ssec:apc-nd}, and \ref{ssec:l1} to compute a coefficient vector $c$. We conduct two kinds of tests to assess accuracy: (a) recovery plots, where the right-hand side vector $b$ in \eqref{eq:gPC_preL1} is generated using an $s$-sparse vector $c^\ast$, and we are interested in the empirical probability that $c = c^\ast$ is computed; (b) accuracy plots, where $b$ is generated from evaluations from a given function $f$, and the error between $f$ and $f_N$ in \eqref{eq:finite-N-expan} is computed. The latter error is computed as a discrete norm over $E = 10000$ samples generated iid from $\omega$.

We are interested primarily in investigating how the sampling number $M$ affects recovery and approximation accuracy. Due to the probabilistic nature of the random sampling method, all reported results are averaged over 100 independent tests to reduce the statistical oscillations. In all our figures and numerical tests, we use ``induced distribution" sampling to stand for the new proposed induced measure sampling, ``CSA" to denote the $\ell^1$ minimization approach with equilibrium measure sampling method from \cite{Jakeman_2016generalizedsample}, and ``MC" to denote a non-preconditioned $\ell^1$ optimization approach \eqref{eq:gPC_L1} where the samples are generated iid from $\omega$.

\subsection{Exact sparse polynomial functions coefficients recovery}\label{ssec:recovery}
\begin{figure}[!htbp]
\begin{center}
\centerline{\includegraphics[width=5cm]{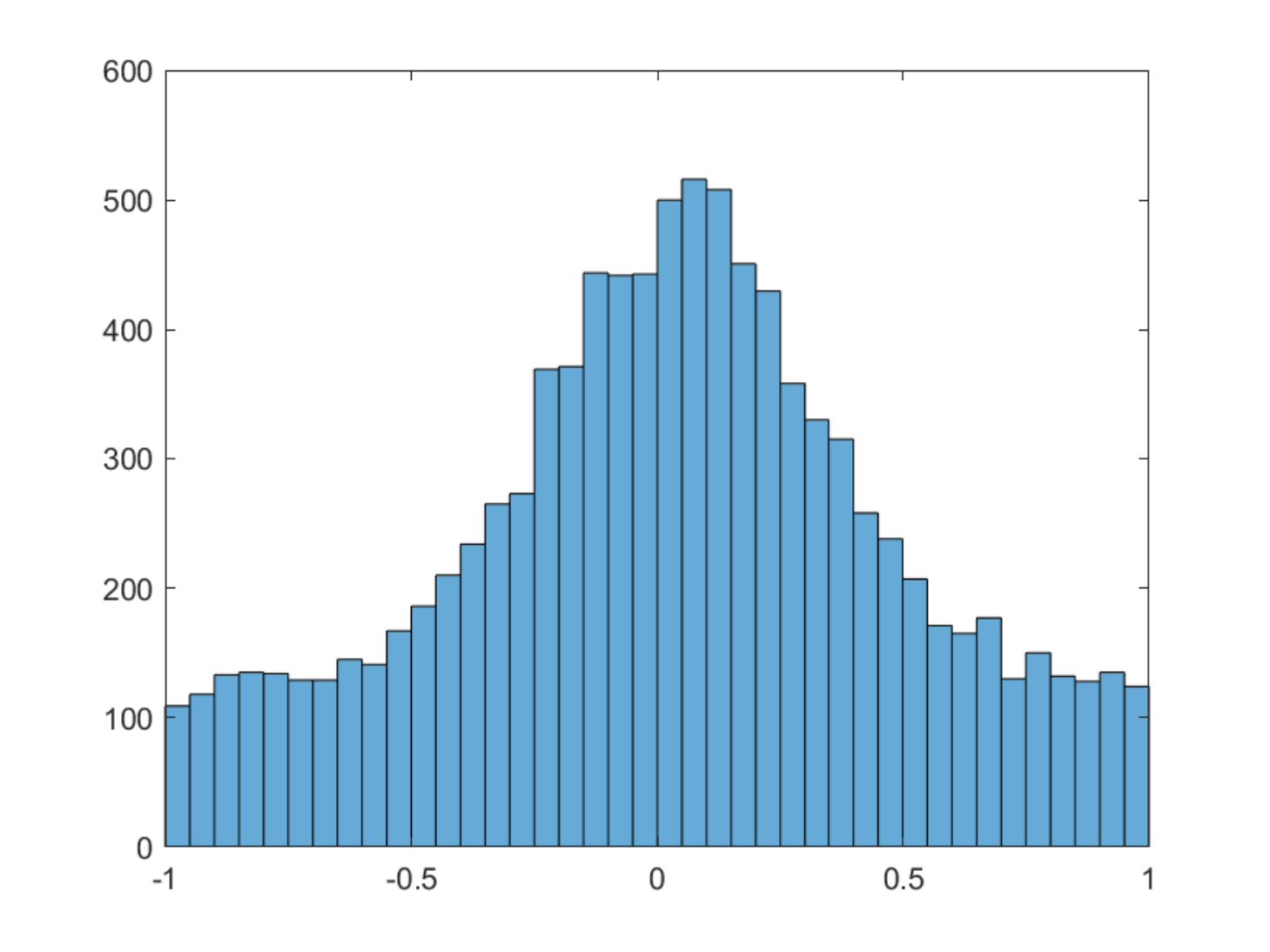}\includegraphics[width=5.0cm]{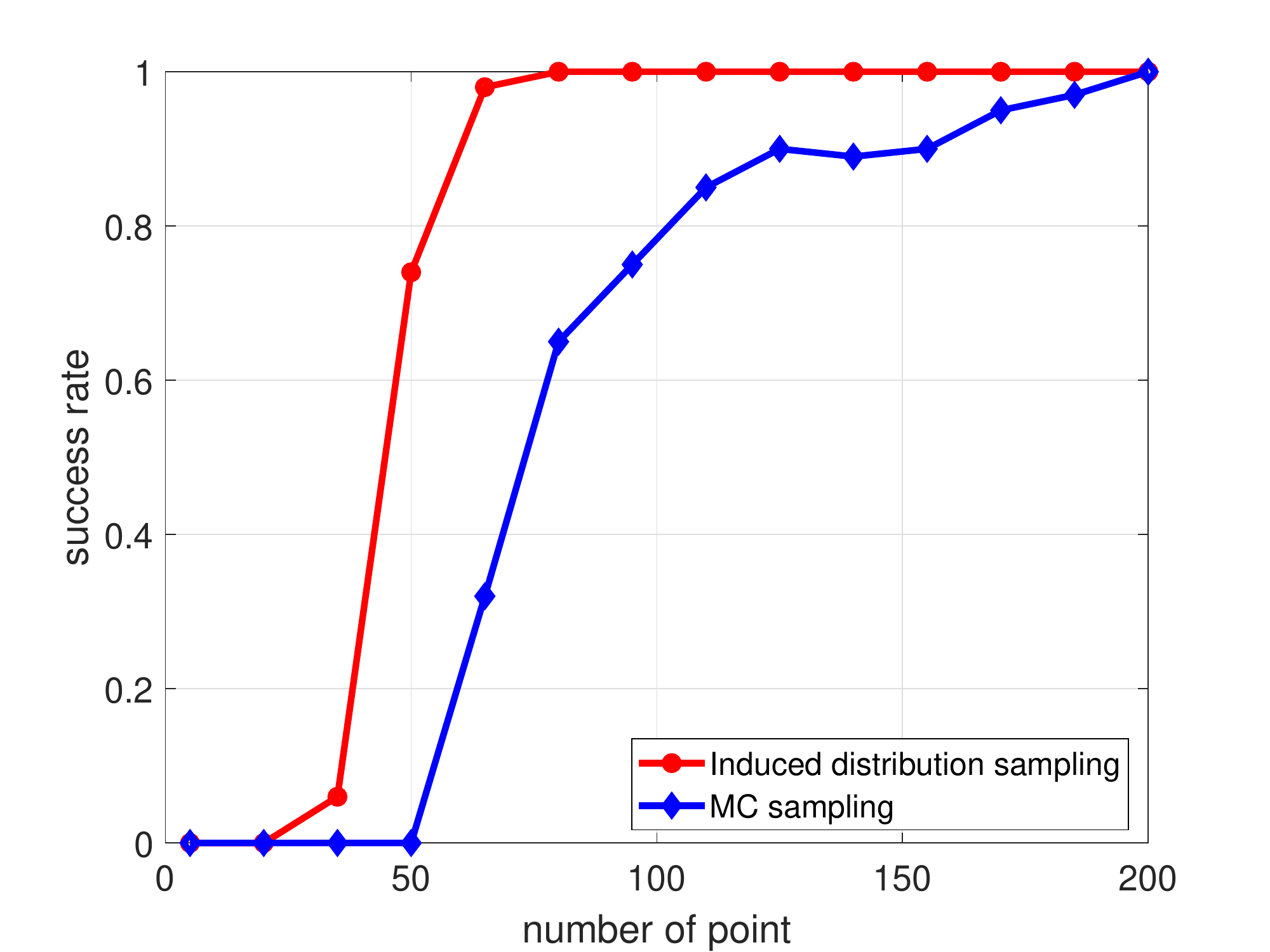}\includegraphics[width=5.0cm]{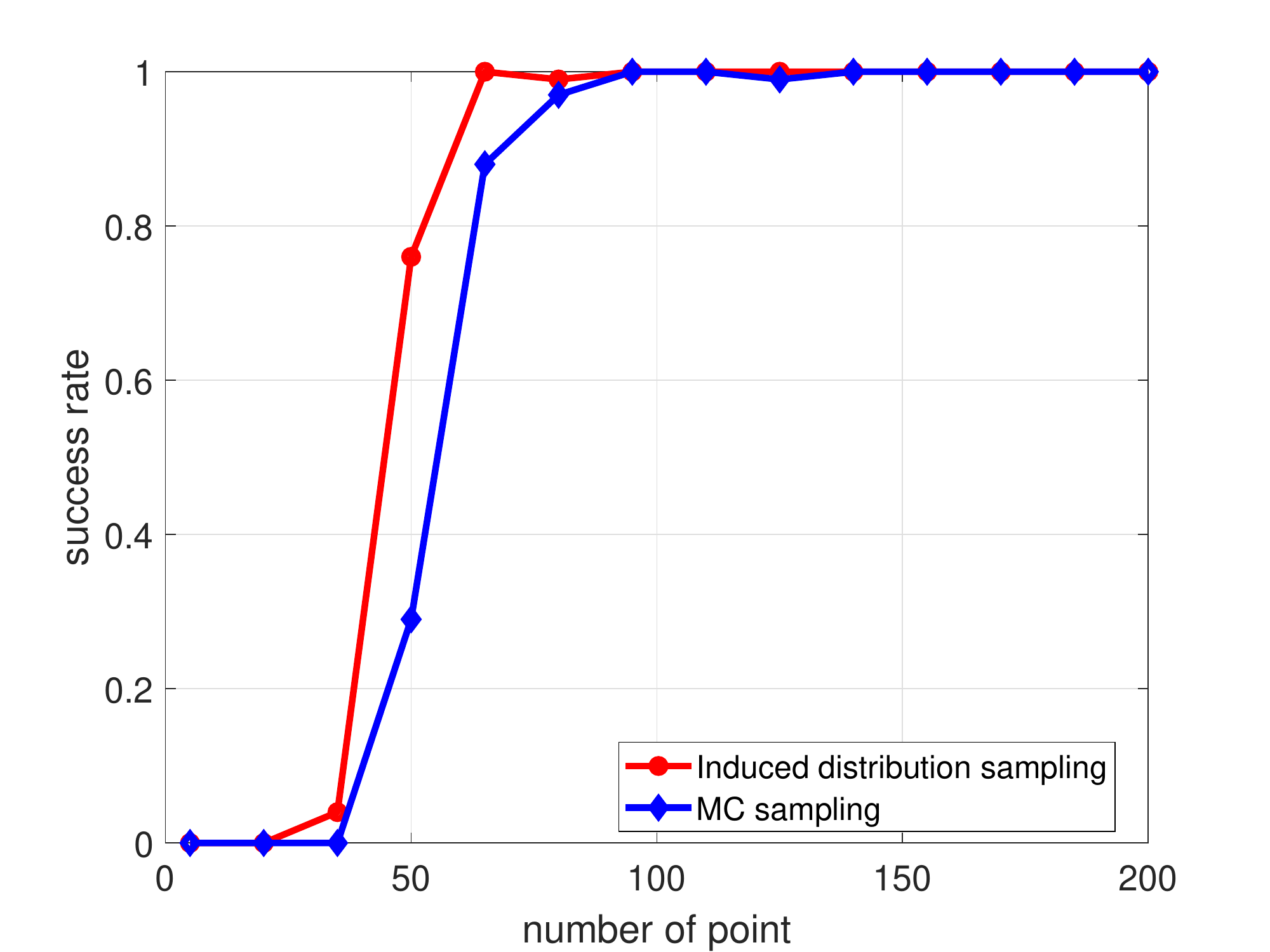}}
\end{center}
  \caption{\sf Left: A histogram of the marginalization of the data $S$. Middle: Probability of successful recovery vs. number of samples $M$ for $d=2$, $K = 20$, $N=231$ with sparsity $s = 8$; Right: Probability of successful recovery vs. number of samples for $d=10$, $K = 3$, $N=286$ with sparsity $s = 8$.}\label{rawdata}
\end{figure}

We first investigate the the performance of the induced sampling method when used to recover manufactured sparse gPC expansions. We set $\omega$ to be a tensorial $d$-dimensional density with isotropic marginal densities. The marginal densities are (equal) mixtures of a density that is uniform on $[-1,1]$, that is a truncated normal density $\mathcal{N}(0.2,1.5)$ on $[-1,1]$, and that is a truncated lognormal distribution on $[0,1]$. We set $S$ to be $Q = 10^5$ iid samples from $\omega$; a one-dimensional histogram for this measure is shown in the left plot of Figure \ref{rawdata}.

Fixing a sparsity level $s$, we choose $c^{\ast}$ to be a randomly generated $s$-sparse vector with standard normal entries. This $s$-sparse vector $c^{\ast}$ is set to be the coefficients of the target arbitrary polynomial expansion, i.e.,
\begin{equation*}
f(z) = \sum_{i=1}^N c^{\ast}_i\Phi_i(z),
\end{equation*}
where $\Phi_i(z)$ is the set of multivariate aPC polynomials generated by $S$. This defines the vector $b$ used in the $\ell^1$ optimization methods. We seek to recover the sparse coefficients via preconditioned $\ell^1$ minimization. In the following experiments, a recovery is considered successful if the resulting coefficient vector $c$ satisfies $\|c-c^{\ast}\|_{\infty}<10^{-3}$.

We first examine the performance of the methods in a relatively low dimension of $d = 2$, the polynomial space is fixed at order $K = 20$ and the number of total unknowns $N=231$. The probability of successful recovery with respect to the number of samples $M$ is displayed in Figure \ref{rawdata} (Middle), for a fixed sparsity level $s = 8$. Results for a higher dimensional case, $d = 10$, with polynomial degree $K = 3$ and a total number of total unknowns $N=286$ is shown in the right plot of Fig. \ref{rawdata}. For the low-dimensional and relatively high-degree situation considered, the induced sampling has a high rate of recovery and performs significantly better than than MC sampling. For the high dimensional and low order case, the induced sampling approach shows similar performance with MC sampling, which is due to low degree polynomials forming an induced distribution $\mu$ that is not too different from $\omega$.

\subsection{Analytical functions approximations}
In this Section, we test the approximation accuracy of the proposed approximation strategy by evaluating the discrete $\ell_2$-error between the true $f$ and the constructed $f_N$ to measure the performance of the approximation. Given a function $f(z)$ and a set of random samples $\{z^{(j)}\}_{j=1}^{E}$, we evaluate the numerical error via
\begin{align*}
  \varepsilon =\left(\frac{1}{M}\sum_{j=1}^{E}|f_N(z^{(j)})-f(z^{(j)})|^2\right)^{1/2},
\end{align*}
where $f_N$ is the approximation \eqref{eq:finite-N-expan} obtained by the data-driven preconditioned $\ell^1$ approach \eqref{eq:gPC_preL1}.
We will consider the following different test functions, which are standard test functions for multivariate approximation \cite{genz_testing_1984}:
\begin{align*}
&f_1(z)=\exp\left(-\sum_{i=1}^dz_i\right),\quad f_2(z)=\sum_{i=1}^{d}(1-z_{i-1})^2+\sum_{i=2}^{d}100(z_{i}-z_{i-1}^2)^2,\\
&f_3(z)=\sin\left(\sum_{i=1}^dz_k\right), \quad f_3(z)=\left(1+\frac{1}{2d}\sum_{i=1}^dc_i(1+z_i)\right)^{-d-1}, \ c_i=\frac{1+i}{4d}.
\end{align*}

\begin{figure}[ht!]
\centerline{\includegraphics[width=7.0cm]{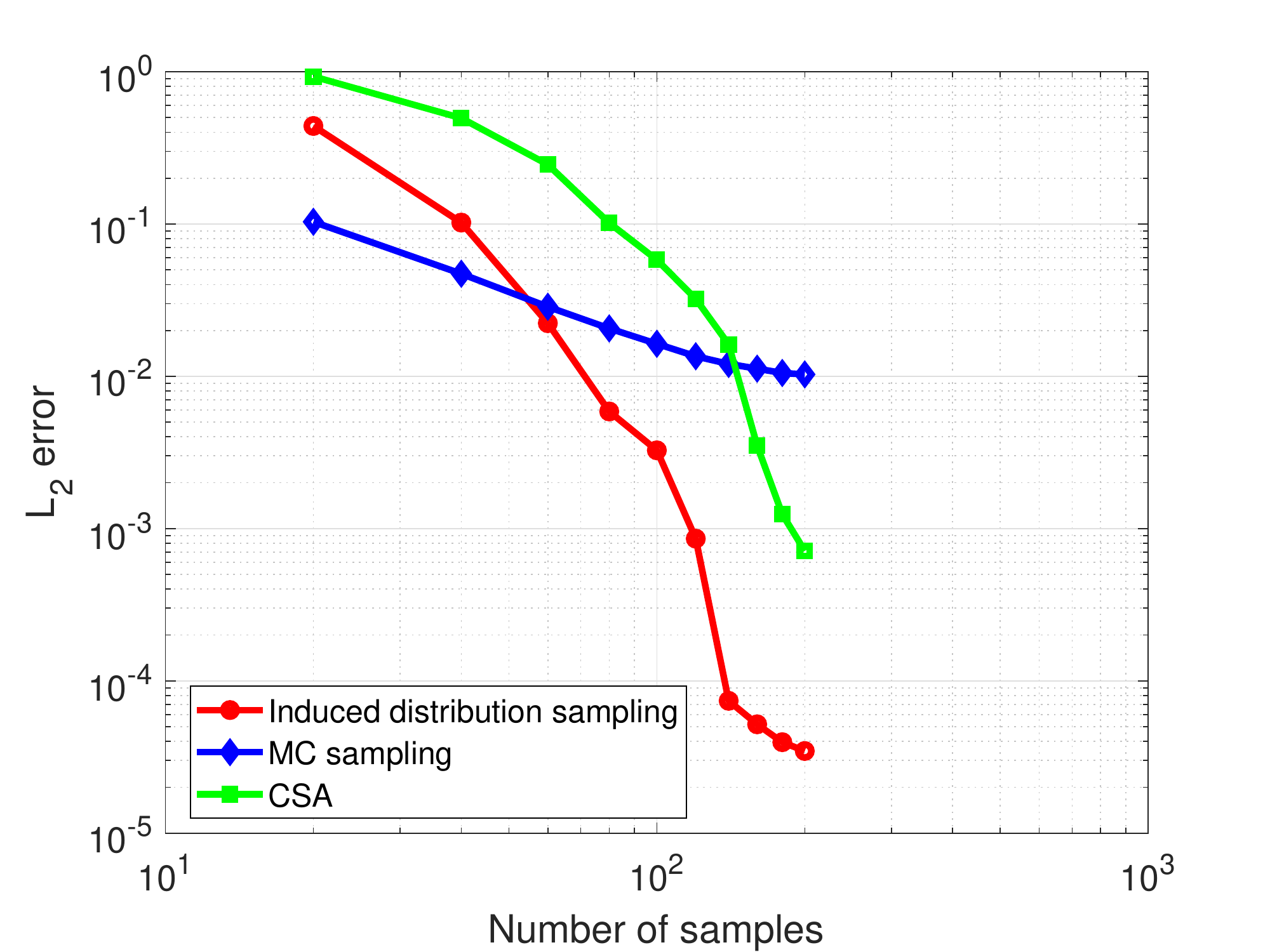}\includegraphics[width=7.0cm]{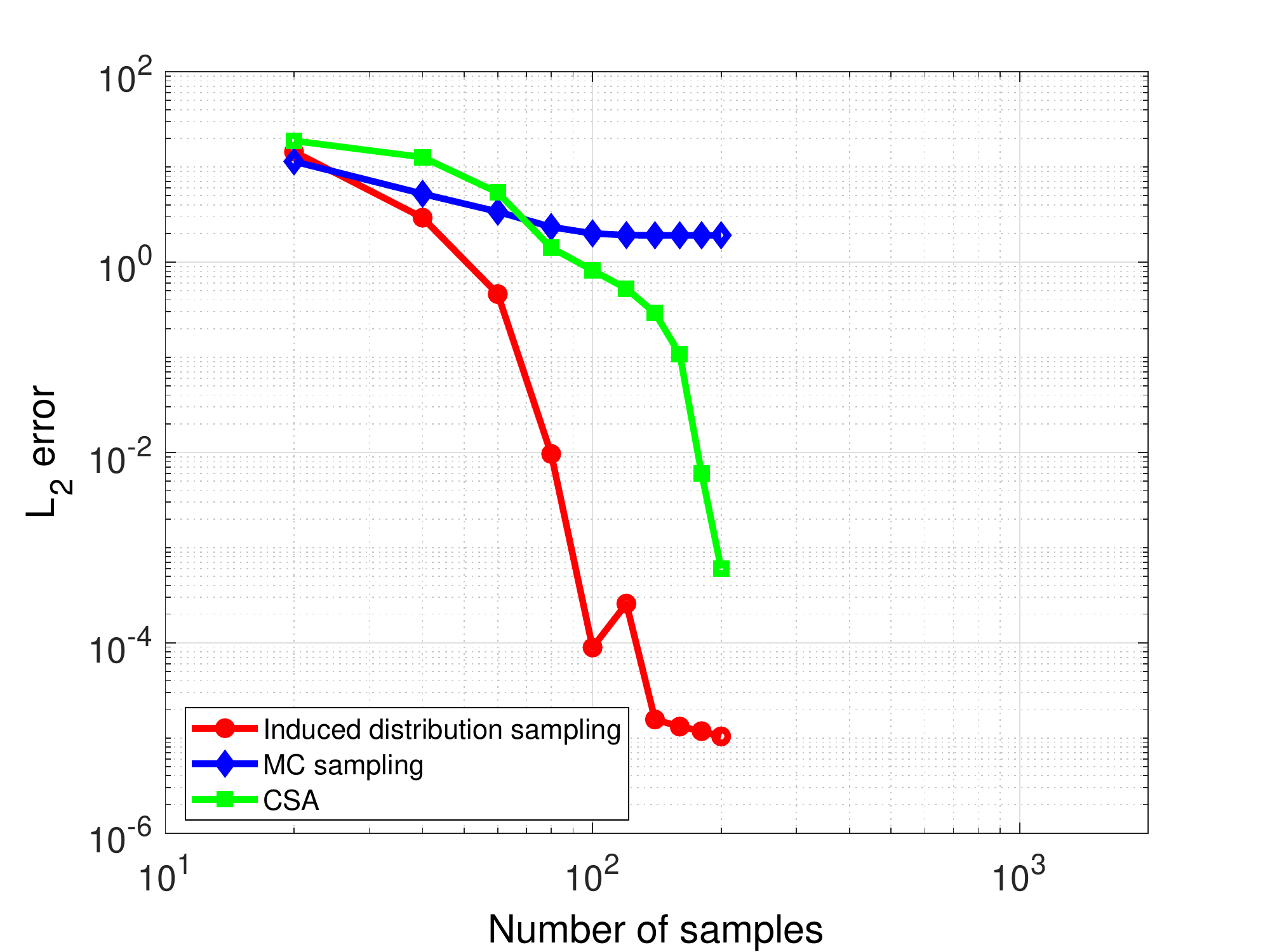}}
\centerline{\includegraphics[width=7.0cm]{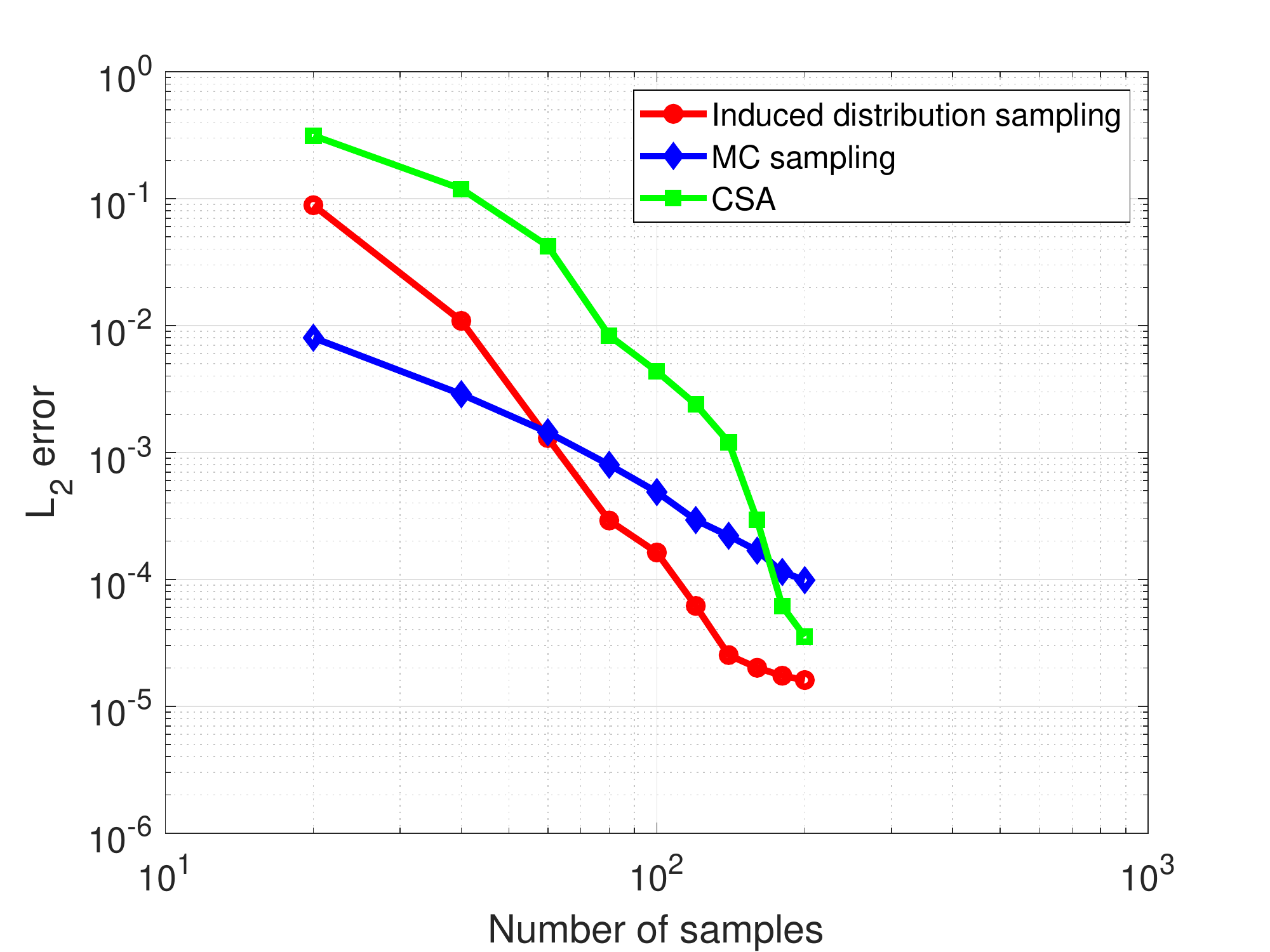}\includegraphics[width=7.0cm]{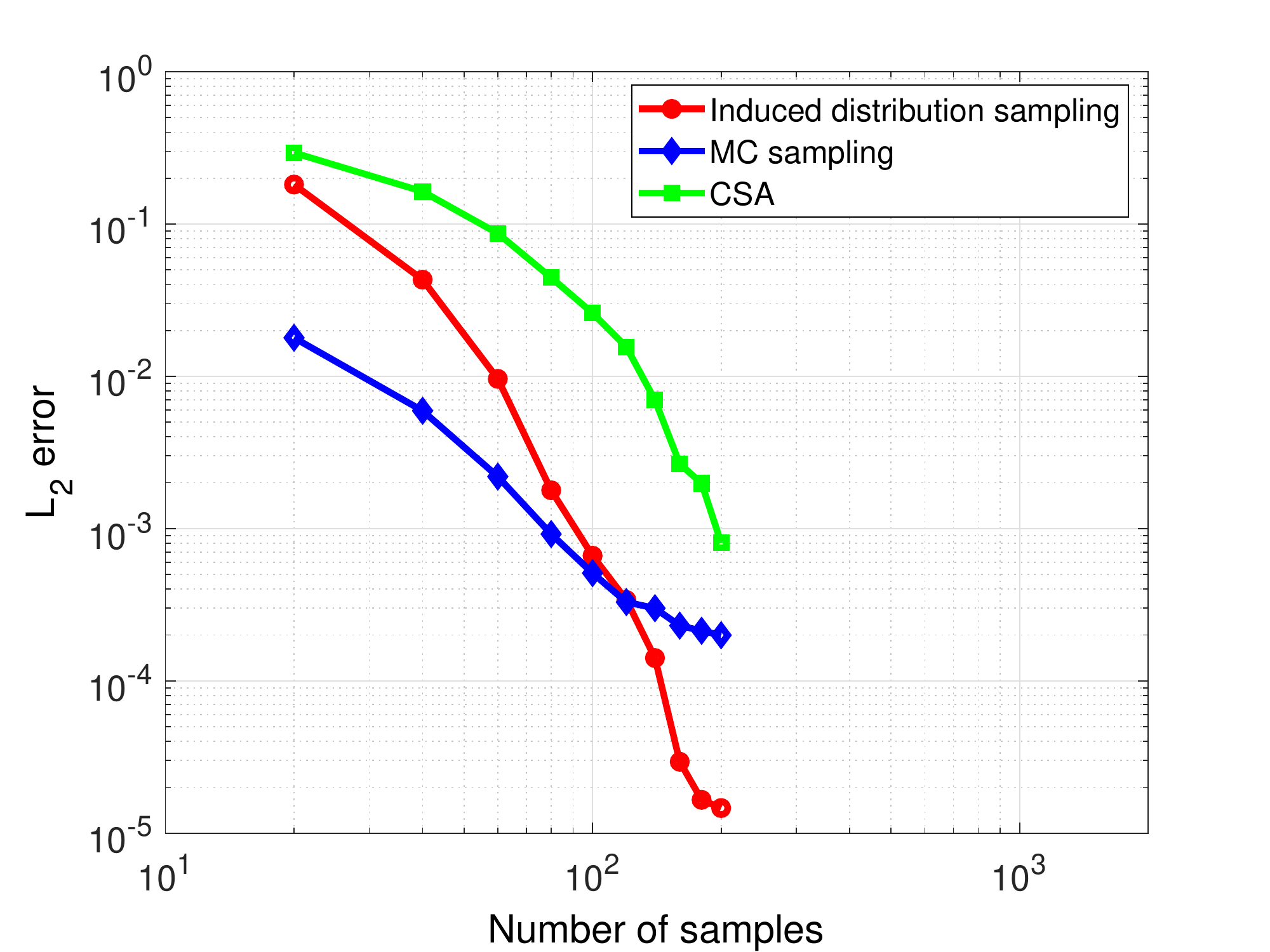}}
  \caption{\sf The discrete $\ell_2$ error of function approximation versus $M$ for three different kinds of sampling strategy. The data set emulates a tensorial measure $(z_1,z_2)\sim \mathrm{Binomial}(24,0.5)\bigotimes \mathrm{Poisson}(10)$. Top Left: $f_1(z)$; Top Right: $f_2(z)$; Bottom Left: $f_3(z)$; Bottom Right: $f_4(z)$.}\label{func2dBP}
\end{figure}

\subsubsection{Two-dimensional case}
In this $d=2$-dimensional case, the sample set $S$ is as constructed in Figure \ref{sample comparison} and \eqref{eq:S-example}.  We set the polynomial space order as $K=20$. 
We use this example to illustrate that the equilibrium sampling approach is not a good choice for this case since it generates samples corresponding to locations of low probability, cf. Figure \ref{sample comparison}. In Figure \ref{func2dBP}, we plot the discrete $\ell_2$ error computed by the $\ell^1$ optimization approach based on induced aPC, CSA, and MC sampling strategies. We observe that, for all test functions, the induced distribution sampling method obtains much better approximation results than that of the MC and CSA sampling approaches for larger values of $M$. For relatively small values of $M$, the the MC method moderately outperforms the aPC induced approach. 

\begin{figure}[ht!]
\centerline{\includegraphics[width=7.0cm]{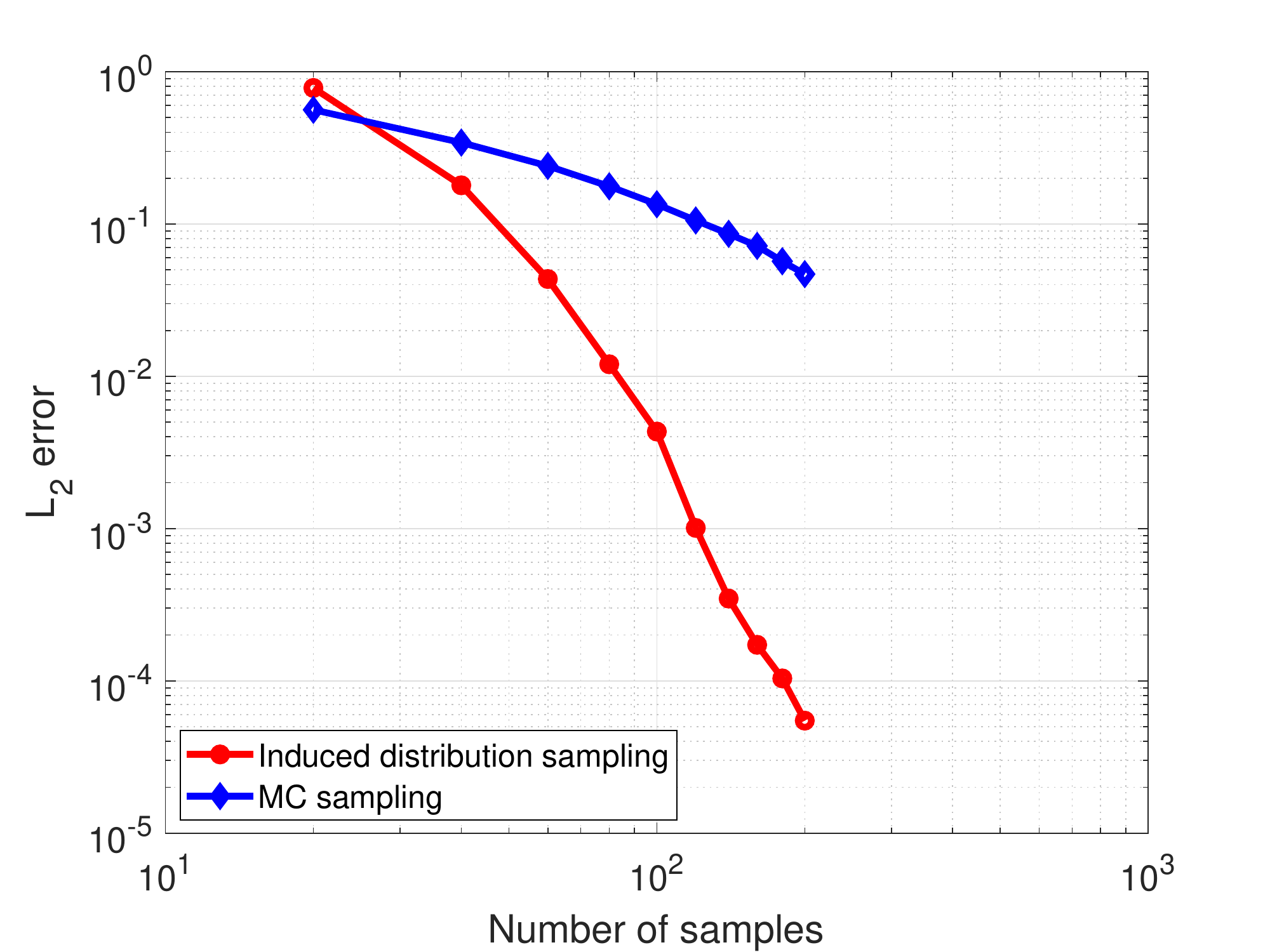}\includegraphics[width=7.0cm]{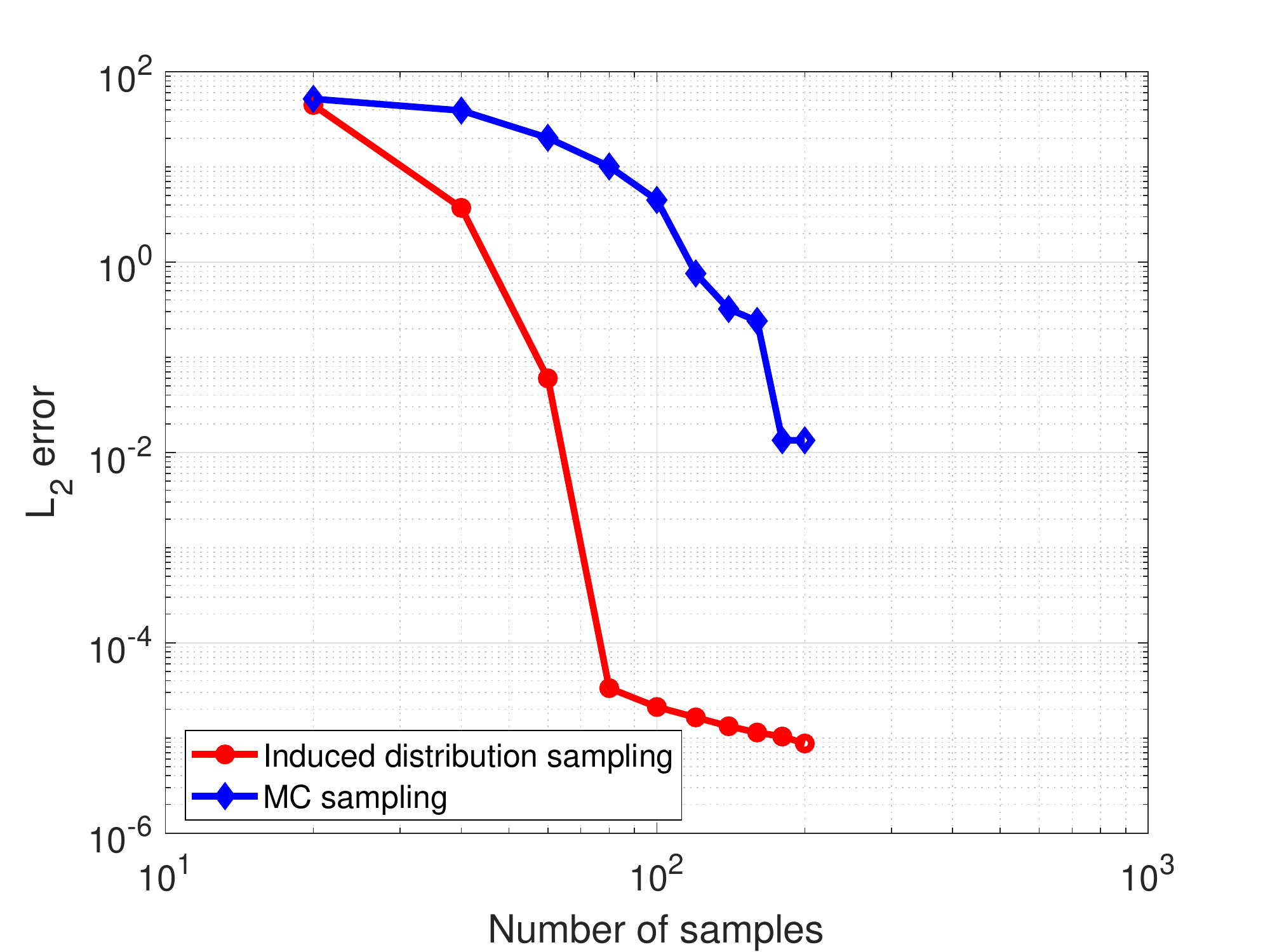}}
\centerline{\includegraphics[width=7.0cm]{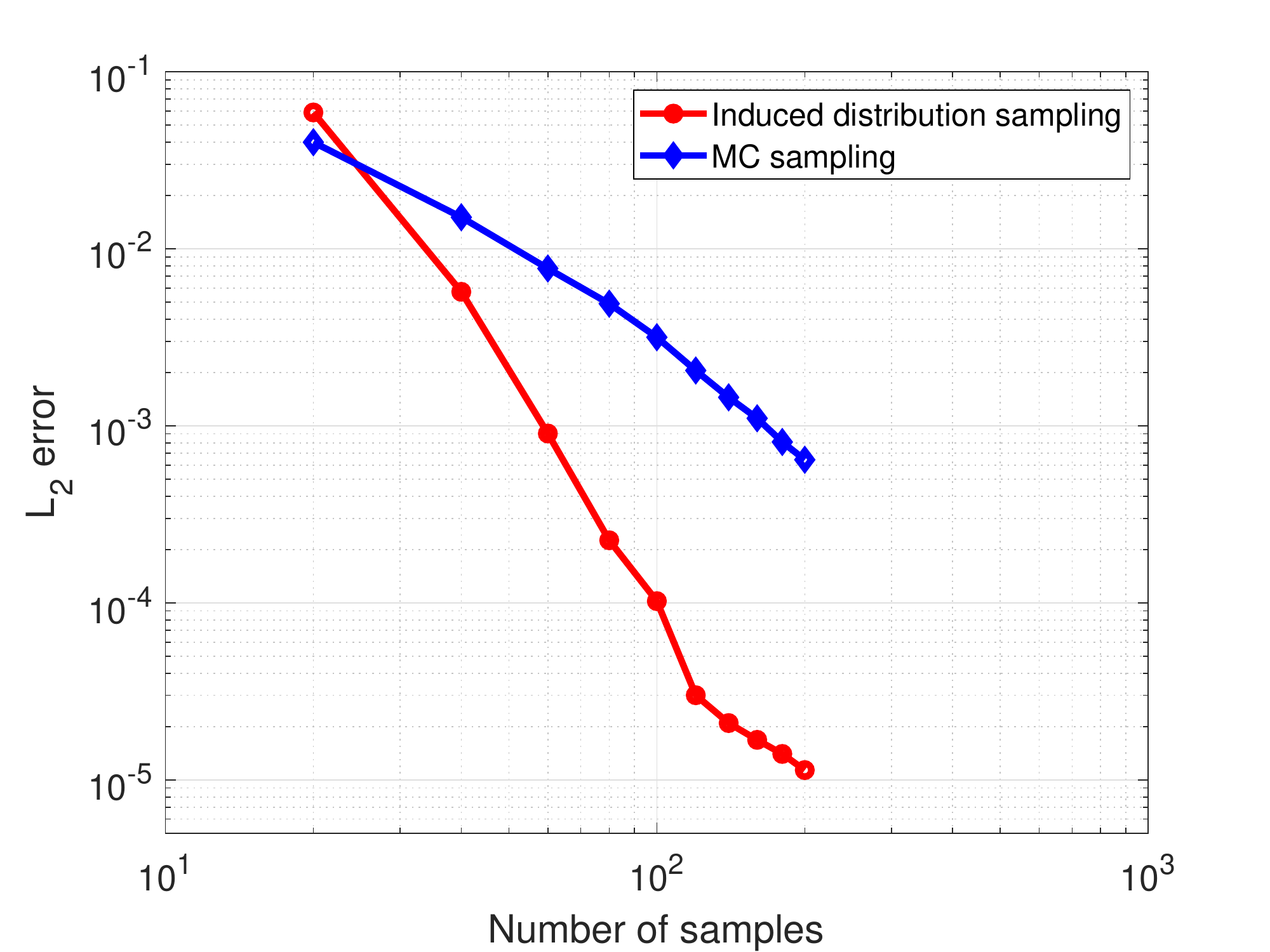}\includegraphics[width=7.0cm]{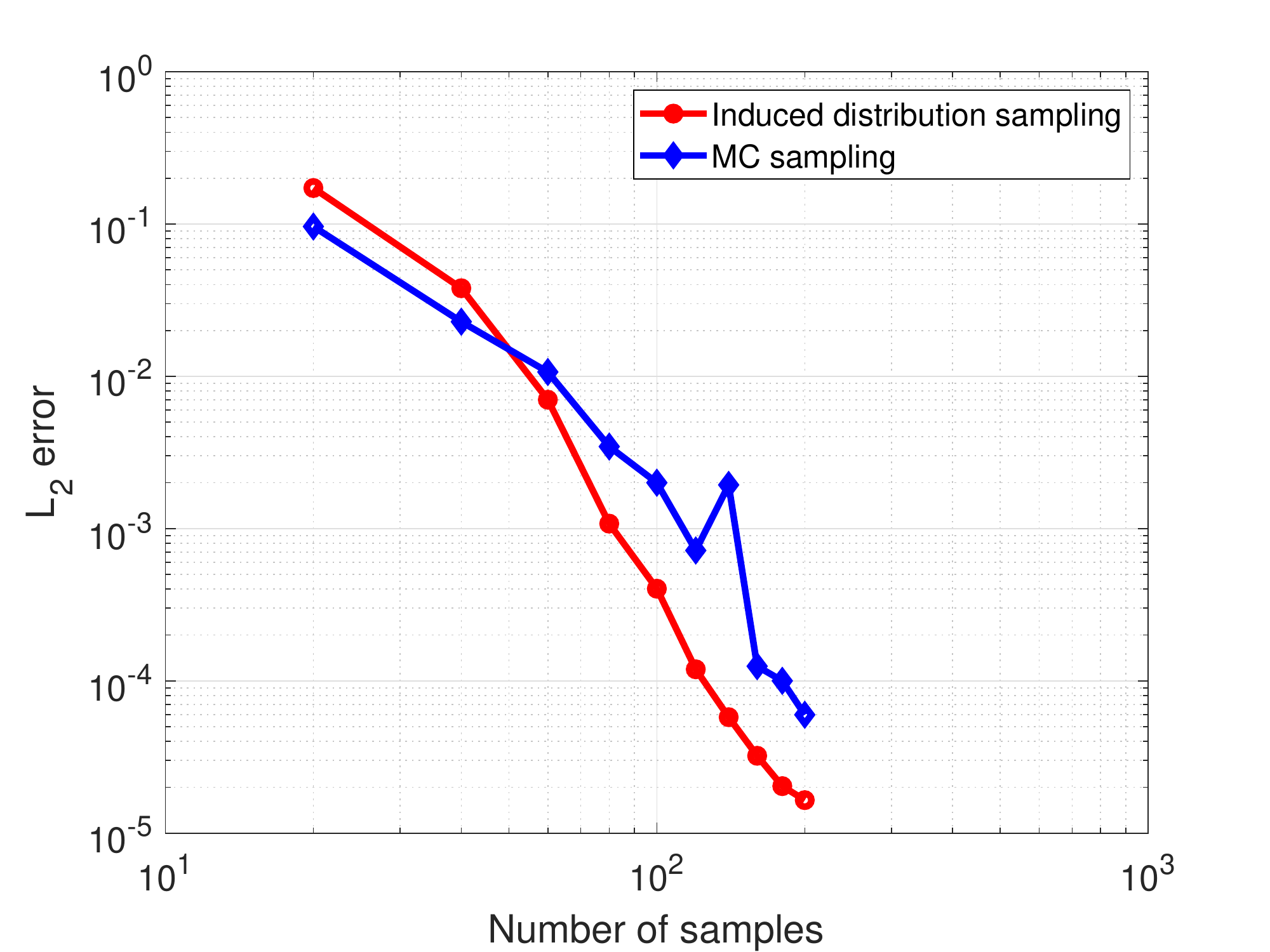}}
\caption{\sf Comparison of the discrete $\ell_2$ error of two-dimensional function approximation versus $M$. The data set is iid from a two-dimensional dimensionally isotropic product density. Top Left: $f_1(z)$; Top Right: $f_2(z)$; Bottom Left: $f_3(z)$; Bottom Right: $f_4(z)$.}\label{func2draw}
\end{figure}

To further study the performance of the proposed induced sampling strategy in the data-driven case, we change $S$ and $\omega$ to correspond to the experiment in Section \ref{ssec:recovery} as shown in the histogram in Figure \ref{rawdata}, left. We show the analytical function approximation errors for this case in Figure \ref{func2draw}. Again, the proposed aPC induced distribution sampling approach obtains much more accurate results than the MC method for increasing $M$.

\subsubsection{Five-dimensional case}
\begin{figure}[ht!]
\centerline{\includegraphics[width=7.0cm]{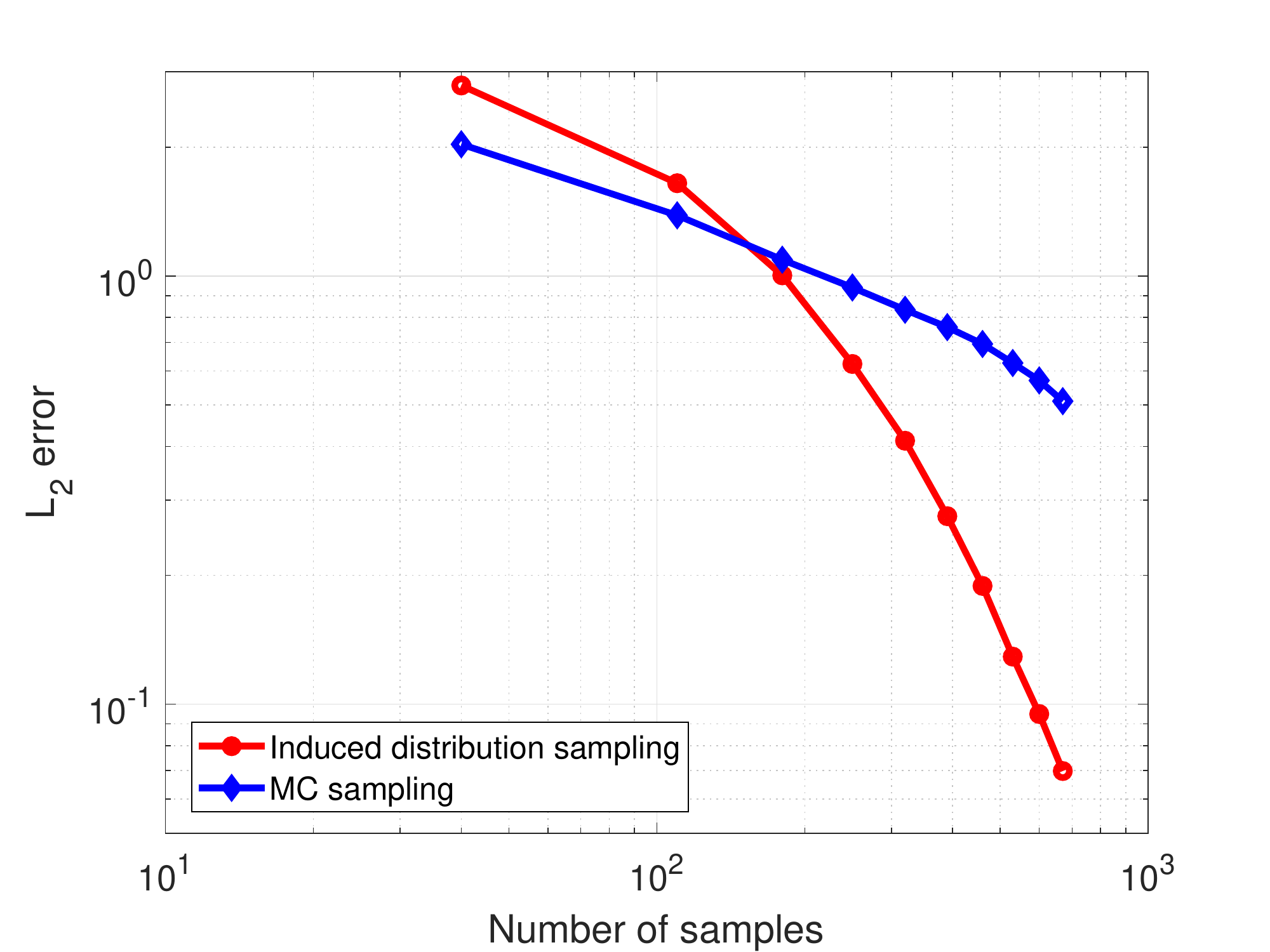}\includegraphics[width=7.0cm]{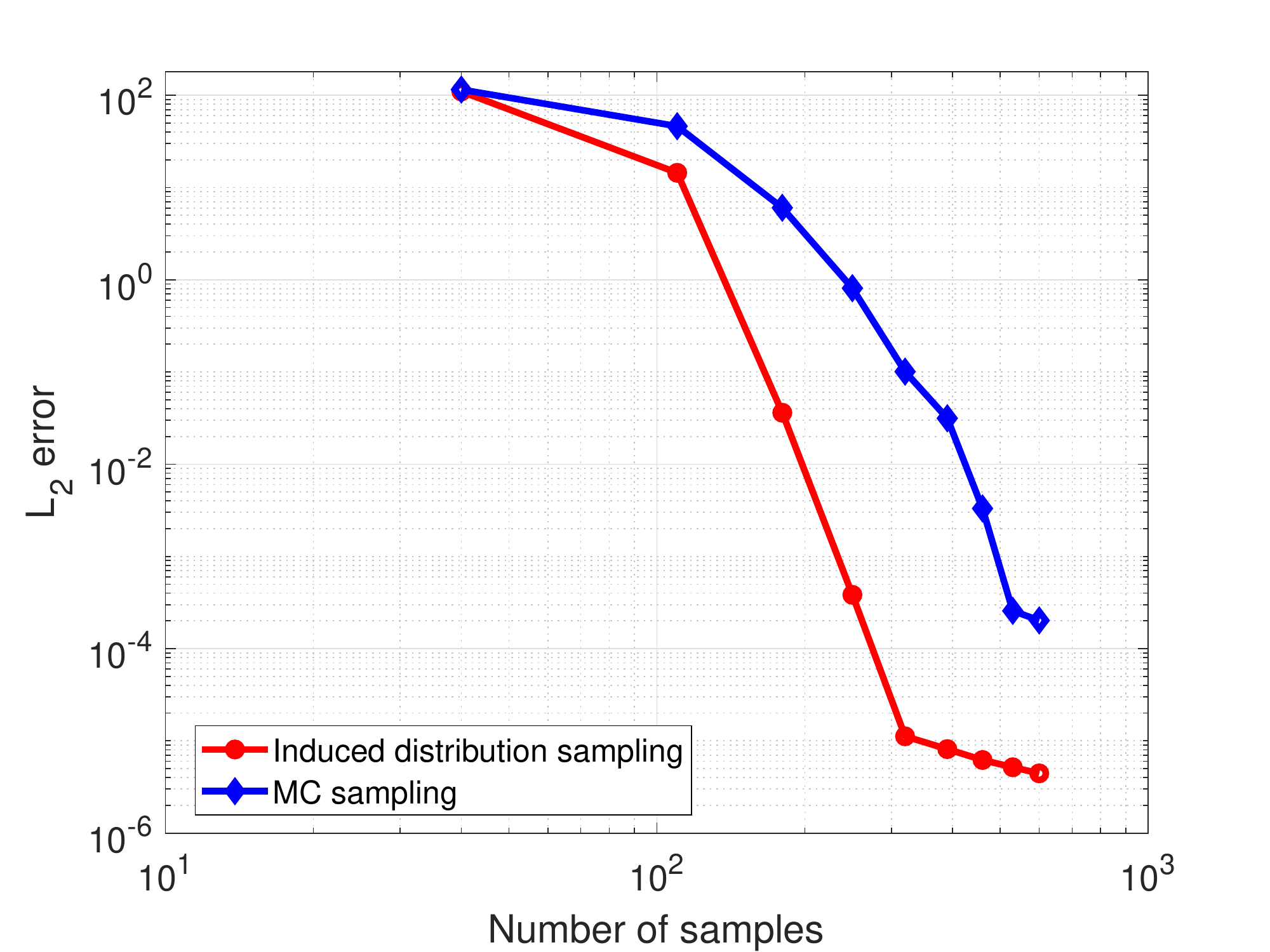}}
\centerline{\includegraphics[width=7.0cm]{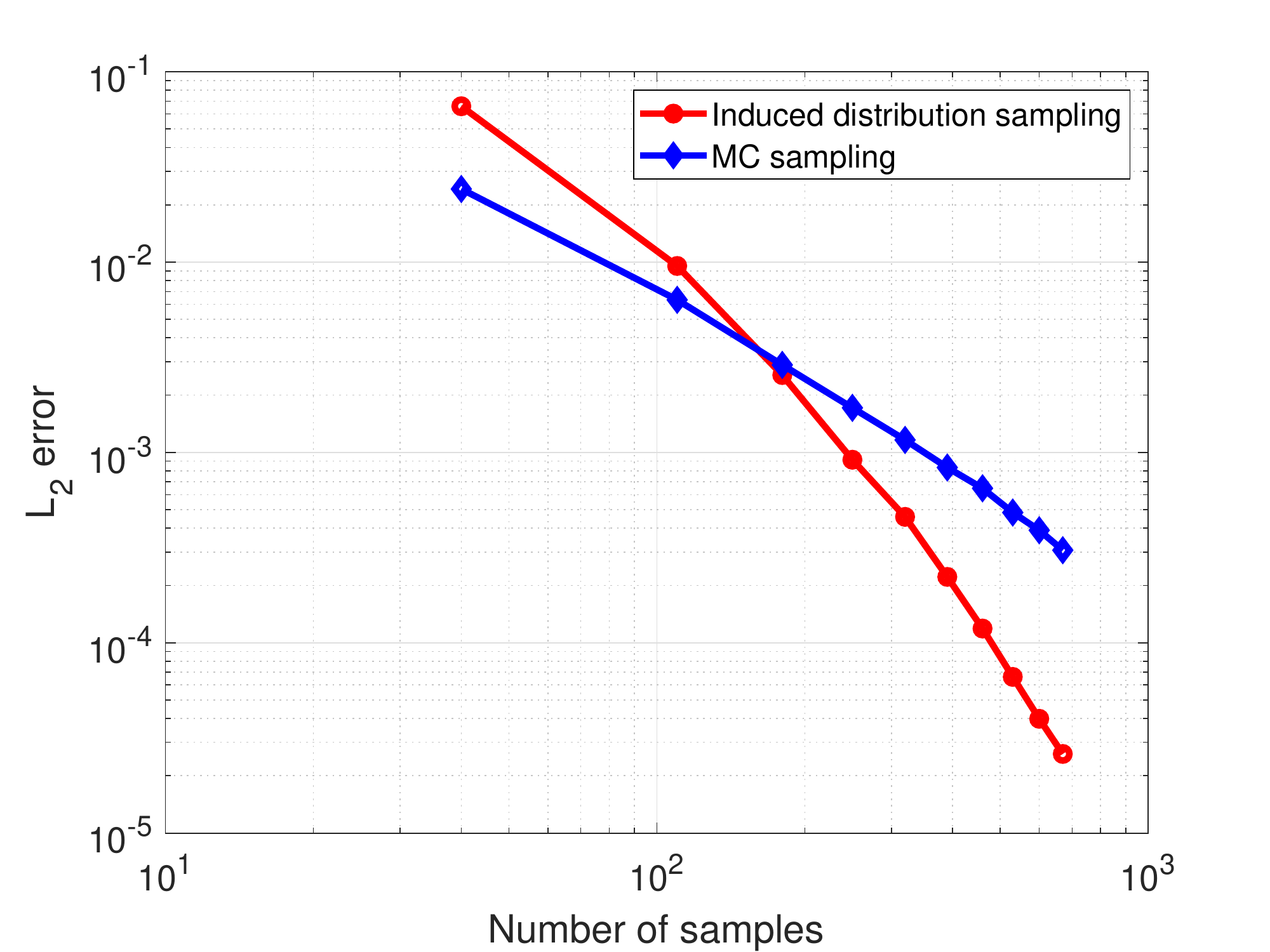}\includegraphics[width=7.0cm]{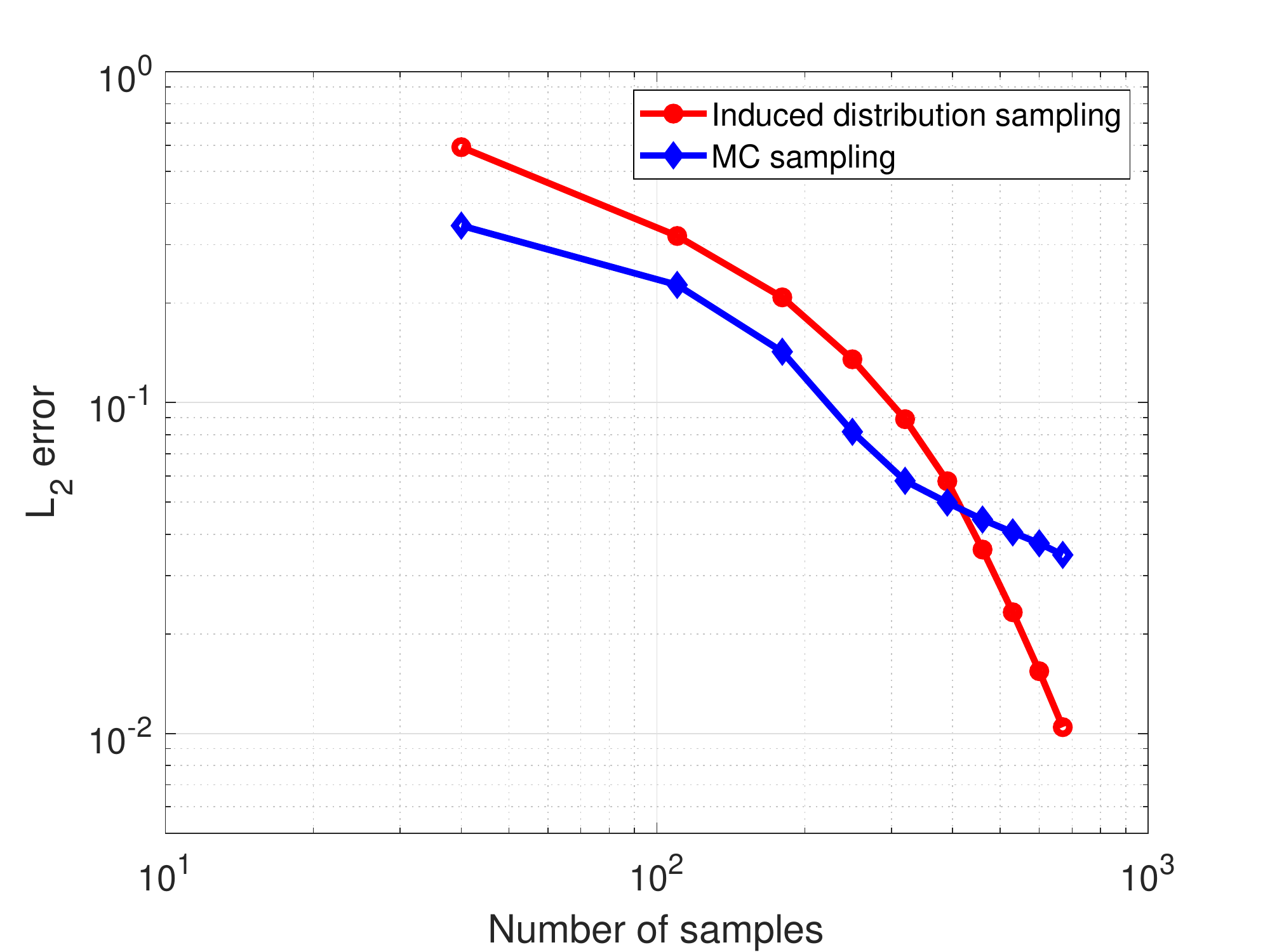}}
\caption{\sf Comparison of the discrete $\ell_2$ error of five-dimensional function approximation versus $M$. The data set is iid from a five-dimensional dimensionally isotropic product density. Top Left: $f_1(z)$; Top Right: $f_2(z)$; Bottom Left: $f_3(z)$; Bottom Right: $f_4(z)$.}\label{func5draw}
\end{figure}

In the five dimensional case, we set $S$ and $\omega$ as in Section \ref{ssec:recovery} and Figure \ref{rawdata}, left.  The polynomial space order is $K = 7$. In Figure \ref{func5draw}, approximation results are shown. We again conclude that the induced distribution sampling method performs better than the MC method for increasing $M$ even in this five dimensional case.

\subsection{Plate bending with random inputs}
Our final example considers the approximation of parameterized partial differential equations, with our goal approximation with respect to the parameter.
Let $\mathscr{D}$ be a bounded polygonal domain in $\mathbb{R}^2$,
with points having the Cartersian coordinate representation
$\boldsymbol{x}=\{x_1,x_2\}^T$. We consider the following clamped Kirchhoff plate bending
problem, which depends on random parameters $z$,
\begin{equation}
\left\{
\begin{array}
{ll}%
-\mathcal{M}_{IJ,IJ}(u(z,\boldsymbol{x}))=f(z,\boldsymbol{x}) &\text{in} \ \ \ \Gamma\times\mathscr{D},\\
u(z,\boldsymbol{x})=\partial_{\boldsymbol{n}}u(z,\boldsymbol{x})=0 & \text{on} \ \ \ \Gamma\times\partial \mathscr{D},\\
\end{array}
\right.  \label{sde}%
\end{equation}
where $\boldsymbol{n}$ denotes the unit outward normal to $\partial\mathscr{D}$, and
\begin{align*}
&{{\mathcal{M}}}_{IJ}(u(z,\boldsymbol{x})):=D(z,\boldsymbol{x})((1-\nu){{\mathcal{K}}}_{IJ}(u(z,\boldsymbol{x})
)+\nu{{\mathcal{K}}}_{LL}(u(z,\boldsymbol{x}))\delta_{IJ}),\nonumber\\
& {{\mathcal{K}}}_{IJ}(u(z,\boldsymbol{x})):=-\partial_{IJ}u(z,\boldsymbol{x})%
=-\frac{\displaystyle
\partial^{2}u(z,\boldsymbol{x})}{\displaystyle \partial
x_{I}\partial x_{J}},I,J=1,2,
\end{align*}
with $D(z,x)=\frac{E(z,\boldsymbol{x})h^3}{12(1-\nu^2)}$ being the rigid flexibility of the plate, $E(z, \bs{x})$ denotes the Young's modulus, $\nu$ is Poisson's ratio, $h$ stands for the thickness of the plate, $\delta_{IJ}$ is Kronecker delta function, and $f(z,x)$ denotes the load force. Equation \eqref{sde} uses the Einstein summation convention, where repeated indices $I$, $J$, and/or $L$ are summed over.

We assume that the uncertainty is imposed on the Young's modulus and satisfies $Y(z, \bs{x})=\log(E(z, \bs{x})-100)$, where $Y(z, \bs{x})$ is given in terms of eigenfunctions for the squared exponential covariance kernel along $x_1$ direction,
$$K(x_1, x'_1)=\exp(\frac{-(x_1-x'_1)^2}{L_c^2}).$$
These eigenfunctions are given by
\begin{equation*}
g_{i}(\bs{x}):= \begin{cases}
\sin\Big(\frac{-(\lfloor\frac{i}{2}\rfloor\pi  x_{1}}{L_{p}}\Big), \ \ i \ \, \textmd{even},\\[12pt]
\cos\Big(\frac{-(\lfloor\frac{i}{2}\rfloor\pi  x_{1}}{L_{p}}\Big), \ \ i \ \, \textmd{odd},
\end{cases}
\end{equation*}
which depend only on $x_1$. Thus $Y(z, \bs{x})$ is one-dimensional
spatial dependent and does not depend on $x_2$.
We will define $Y$ through a truncated Karhunen-Lo\`{e}ve expansion, i.e.,
$$Y = Y(z,\bs{x})=1+Z_1\Big(\frac{\sqrt{\pi}L}{2}\Big)^{1/2}+\sum_{i=2}^{d} \zeta_{i}g_{i}(\bs{x})Z_{i},$$
where
\begin{equation*}
\zeta_{i}:=(\sqrt{\pi}L)^{1/2}\exp\Big(\frac{-(\lfloor\frac{i}{2}\rfloor\pi L)^{2}}{8}\Big),\, \ \textmd{for } \ \,i>1
\end{equation*}
and $\{Z_{i}\}^{d}_{i=1}$ are random variables that we assume are independent. We set $L_{c}=1/2$ and the parameters $L_{p}$ and $L$ are defined as $L_{p}=\max\{1,2L_{c}\}$ and $L=\frac{L_{c}}{L_{p}}$, respectively.

The plate domain is $\mathscr{D}=[0,1]^2$ and it is subjected to a deterministic load $f(z,\boldsymbol{x})=\cos(x_1)\sin(x_2)$, and the plate thickness is set as $h=0.02$. The deterministic plate bending problem is solved by the Morley nonconforming finite element method with the space domain $\mathscr{D}$ been partitioned into 648 triangles corresponding to 1369 unknowns. The quantity of interest that we seek to approximate with respect to $z$ is $u(z) = u(z, (0.5, 0.5))$.

\begin{figure}[htbp]
\begin{center}
    \includegraphics[width=6cm]{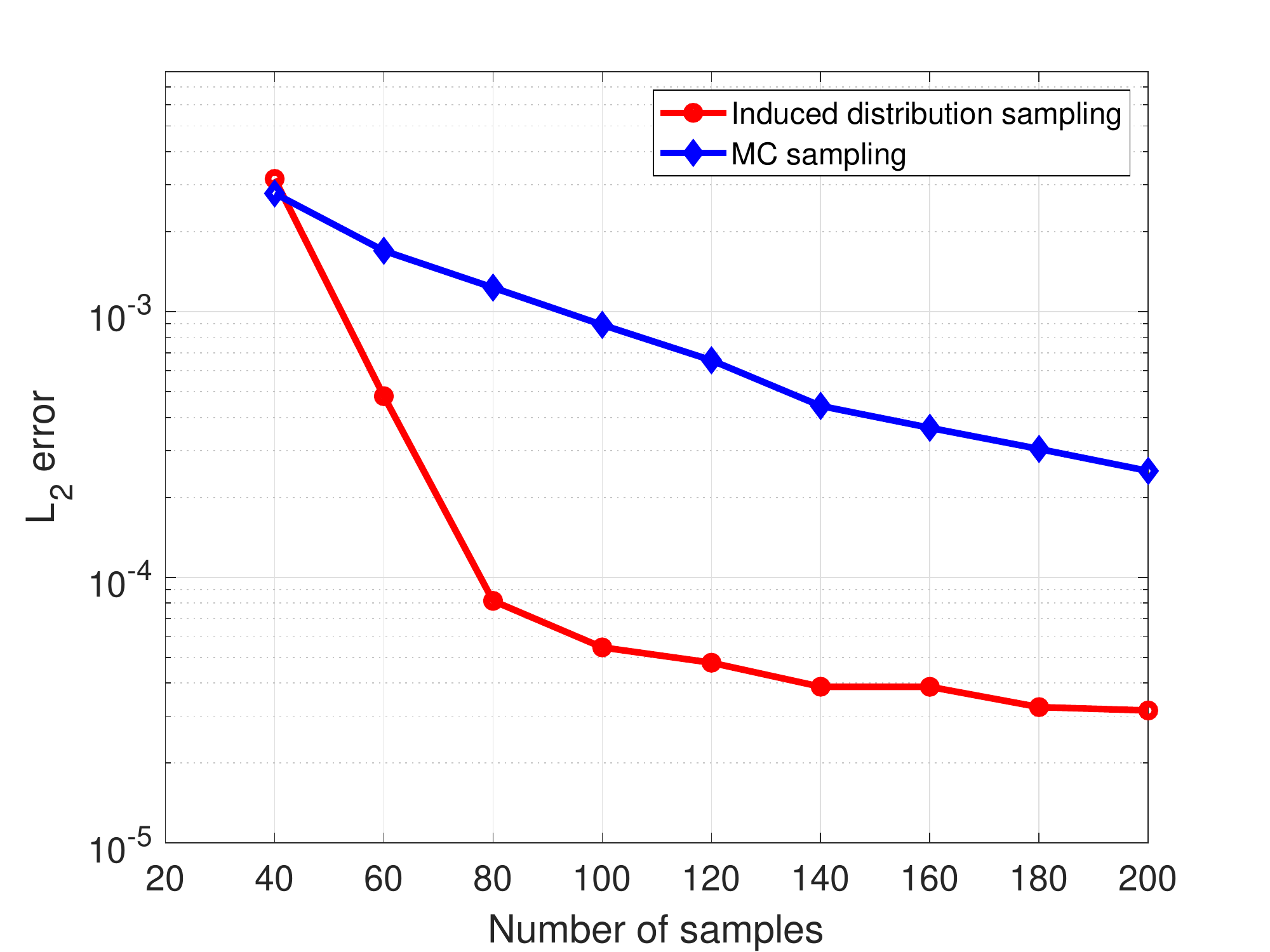}
    \includegraphics[width=6cm]{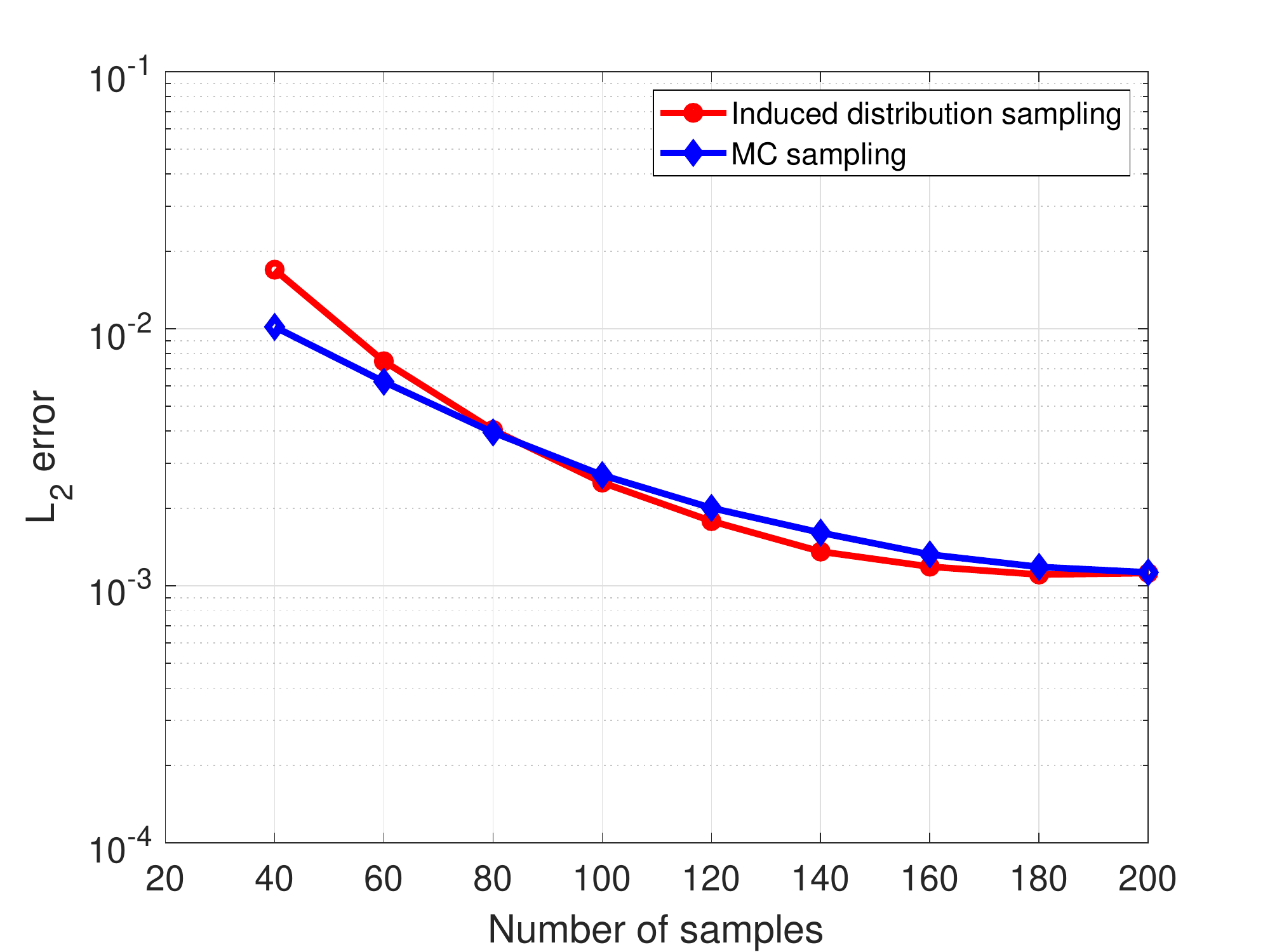}
\end{center}
\caption{The effect of dimension on the performance of the induced sampling aPC approximation
of the plate bending equation \ref{sde}. (Left) 20th degree polynomial in 2 dimensions. (Right) 3th degree
polynomial in 10 dimensions. }\label{fig:PDE_err}
\end{figure}

We consider the two and ten dimensional examples, where we assume that $\omega$ is tensorial and dimensionally isotropic with marginal density as in Section \ref{ssec:recovery} and Figure \ref{sample comparison}, left. We use $Q = 10^4$ iid samples from $\omega$ to define $S$. Approximation results are shown in Figure \ref{fig:PDE_err}. For the low-dimensional case ($d=2$), the provided induced sampling significantly out-performs Monte Carlo sampling: Many fewer samples are needed with the aPC induced approach to achieve a fixed $l_2$ accuracy. However, for the high-dimensional case ($d = 10$), MC achieves accuracy comparable to the induced distribution sampling. This is due to the fact in our high-dimensional approximations we can only use a low-degree polynomials, which inhibits the effectiveness of the induced sampling method (since $\mu$ behaves very close to $\omega$ in this case). These findings are consistent with the behavior observed in other studies on sparse approximation for gPC expansions via different sampling method \cite{Yan_2012Sc,Hampton_2015Cs,Jakeman_2016generalizedsample}.

\section{Summary}
In this paper, we present a sparse arbitrary polynomial approximation approach to solve UQ problems with random inputs whose distribution is characterized only by a set of empirical samples. The moment matching method developed in \cite{Ahlfeld_2016SAMBA} is employed to construct the so-called arbitrary polynomial basis functions. We subsequently propose to use induced distribution sampling defined by the aPC basis functions to collect a small amount of data. A preconditioned $\ell^1$ optimization problem is then utilized to compute a sparse polynomial approximation. Similar to related literature, we observe that the aPC induced method is notably superior to alternative approaches in low dimensions when a (relatively) high polynomial degree is used.  However, the accuracy effectiveness deteriorates for high dimensional approximation with small polynomial degree. 

\bibliographystyle{plain}
\bibliography{ref}

\end{document}